\documentclass[final,reqno]{amsart}
\usepackage{pdfsync}
\usepackage{amsmath,amssymb}
\usepackage{enumerate}
\usepackage{xspace}
\usepackage{array}
\usepackage{tabularx}
\usepackage{subfigure}
\usepackage{caption}
\usepackage{fullpage}
\usepackage{cite}
\usepackage{caption}
\usepackage{james}

\hypersetup{colorlinks, linkcolor=blue, citecolor=blue, urlcolor=blue}

\usepackage{mathtools}
\mathtoolsset{showonlyrefs=true}

\author{
 Alex Bihlo
}
\address{
  Alex Bihlo
  \thanks{
    Department of Mathematics and Statistics, Memorial University of
    Newfoundland, St.\ John's, NL, A1C 5S7, Canada
    {\tt{abihlo@mun.ca}}.
  }}

\author{
  James Jackaman
}
\address{
  James Jackaman
  \thanks{
    Department of Mathematics and Statistics, Memorial University of
    Newfoundland, St.\ John's, NL, A1C 5S7, Canada
    {\tt{jjackaman@mun.ca}}.
  }}

\author{Francis Valiquette}
\address{
  Francis Valiquette
  \thanks{
    Department of Mathematics, Monmouth University, West Long Branch,
    NJ, 07764, USA
    {\tt{fvalique@monmouth.edu}}.
  }}

\thanks{This research was supported, in part, thanks to the Canada
  Research Chairs, the InnovateNL LeverageR\&D and NSERC Discovery
  grant programs}

\title[Symmetry-preserving finite element methods for ODEs]%
{On the development of symmetry-preserving finite element
  schemes for ordinary differential equations}
\date{\today}
\begin{document}

\begin{abstract}

  In this paper we introduce a procedure, based on the method of
  equivariant moving frames, for formulating continuous Galerkin
  finite element schemes that preserve the Lie point symmetries of
  initial value \revise{problems for} ordinary differential
  equations. Our methodology applies to projectable and
  non-projectable \revise{symmetry group} actions, \revise{to}
  ordinary differential equations of arbitrary order, and
  \revise{finite element approximations} of arbitrary
  \revise{polynomial} degree. Several examples are included to
  illustrate various features of the symmetry-preserving process.  We
  summarise extensive numerical experiments showing that
  symmetry-preserving finite element schemes may provide better long
  term accuracy than their non-invariant counterparts and can be
  implemented on larger elements.
 
\end{abstract}

\maketitle

%%%%%%%%%%%%%%%%%%%%%%%%%%%%%%%%%%%%%%%%%%%%%%%%%%%%%%%%%%%%%%%%%%%%%%%%

\section{Introduction} \label{sec:introduction}

For the accurate long time simulation of ordinary differential
equations (ODEs), it is paramount that intrinsic properties of the
underlying problem be preserved, giving rise to the field of geometric
numerical integration, \cite{HairerLubichWanner:2006}.  Well-known
geometric numerical methods include symplectic integrators
\revise{such as collocation and Runge-Kutta methods},
\cite{BlanesCasas:2016, HairerLubichWanner:2006, LeimkuhlerReich:2004,
  Sanz-SernaCalvo:1994}, Lie--Poisson structure preserving schemes,
\cite{ZhongMarsden:1988}, energy-preserving methods,
\cite{QuispelMcLaren:2008}, general conservative methods,
\cite{WanBihloNave:2016, WanBihloNave:2017}, and more.  By preserving
certain characteristics of the differential equations, these geometric
numerical schemes \revise{typically} provide better global and long term results than
their non-geometric
analogues.

Over the last 30 years there has been a considerable amount of work
focusing on the elaboration of finite difference numerical schemes
that preserve the Lie point symmetries of differential equations,
\cite{BakirovaDorodnitsynKozlov:1997, BuddDorodnitsy:2001,
  Dorodnitsyn:1989, Dorodnitsy:2011, DorodnitsyWinternitz:2000,
  Kim:2008, KimOlver2004}.  Using either infinitesimal generators or
the method of equivariant moving frames, the algorithms for
constructing symmetry-preserving finite difference schemes are now
well-established.  For a review and comparison of the two approaches,
we refer the reader to \cite{BihloValiquette:2017}.  For ordinary
differential equations, symmetry-preserving numerical schemes have
shown to give good numerical results, especially when solutions
exhibit sharp variations or admit singularities, \cite{Bourlioux:2006,
  BourliouxRebeloWinternitz:2008, KimOlver2004}.  For partial
differential equations, the numerical improvements are not as clear,
\cite{BihloCoiteuxWinternitz:2015, Kim:2008,
  LeviMartinaWinternitz:2015, RebeloValiquette:2013}, and there
remains a lot more to be done.  Most of the work thus far has focused
on evolutionary partial differential equations and it is now accepted
that in order to preserve symmetries, numerical schemes have to be
defined on time-evolving meshes.  To avoid mesh tangling and other
numerical instabilities, the basic invariant numerical schemes have to
be combined with evolution--projection techniques, invariant
r-adaptive methods, or invariant meshless discretisations,
\cite{Bihlo:2013, BihloNave:2013, BihloNave:2014, BihloPopovych:2012}.

Recently, the first and third authors began to adapt some of the
techniques used to construct symmetry-preserving finite difference
schemes to finite element methods, \cite{BihloValiquette:2018},
focusing on second order ODEs and Burgers' equation.  From a numerical
perspective, finite element methods offer several advantages over
finite difference methods. \revise{Firstly}, when considering complex domains,
unstructured grids or moving boundaries, finite element methods are
generally easier to implement than finite difference methods.
Secondly, since finite element methods are based on the discretisation
of the weak form of a differential equation, such numerical schemes
have less rigid smoothness requirements than finite difference
schemes.  From an application standpoint, finite element methods are 
used in a wide range of fields such as fluid dynamics, engineering, physics, and
applied mathematics.

Given a group of Lie point transformations acting freely and
regularly on a manifold, the method of equivariant moving frames, originally
developed in \cite{FelsOlver:1999}, and extended to Lie pseudo-groups,
discrete groups, and finite differences in \cite{MariMansfield:2018,
  Olver:2001, Olver:2019, OlverPohjanpelto:2008}, is a powerful tool
for constructing invariant quantities of the group action.  Indeed, the method of equivariant moving frames
is equipped with a systematic process for constructing invariant
quantities using the so-called invariantisation map, which sends
non-invariant quantities to their invariant counterparts.  Therefore, given a finite
element method, which will in general not preserve the symmetries of the original differential equation,
a symmetry-preserving scheme can be easily constructed by invariantising the non-invariant discrete weak
formulation.  This is the same key idea used to construct symmetry-preserving finite difference schemes in
\cite{BihloValiquette:2017, RebeloValiquette:2013, Kim:2008,KimOlver2004},

In \cite{BihloValiquette:2018}, the authors limited themselves to
second order ODEs and only considered projectable group actions. 
In the present paper we generalise the ideas set
out in \cite{BihloValiquette:2018} to ODEs of arbitrary orders and to 
general Lie point transformation groups.  In particular, we utilise the conforming nature of
the finite element approximation to exploit the invariantisation
procedure within a continuous framework, whereas in the prequel a
finite difference like approach was considered. We note that viewing high
degree finite elements in the difference setting leads to 
complex calculations, which adds computational challenges to the
methodology considered in  \cite{BihloValiquette:2018}.  The approach introduced in this paper
has several benefits, not least of which is the ability to construct invariant finite element
approximations of arbitrarily high order simultaneously. Though, we note that
for some group actions it will be necessary to consider certain
aspects of the invariantisation procedure for each polynomial degree
independently.  In these cases, these degree specific
modifications are significantly simpler than those introduced in 
\cite{BihloValiquette:2018}. 

While finite element methods were originally introduced as a
methodology for spatial discretisation, they are also \revise{ideal} for
temporal discretisations. In particular, by rewriting an arbitrarily
high order ODE as a system of first order equations we may introduce
both continuous and discontinuous time-stepping Galerkin
approximations which are well studied, see \cite{Johnson:1988,
  FrenchSchaeffer:1990, EstepFrench:1994, Estep:1995}. It is known
that for a large class of geometric ODEs, the continuous Galerkin
method preserves the underlying energy of the problem,
\cite{Hansbo:2001}, while this is not the case for the discontinuous
Galerkin method, \cite{EstepStuart:2002}.  Therefore, in the sequel we
will disregard the discontinuous method and focus on its conservative
continuous counterpart. \revise{Furthermore, the continuous Galerkin
approximation is known to be ``close'' to a family of symplectic
collocation methods,} see \cite[Chapter 3.1]{self:thesis}.

The remainder of this work is set out as follows: In Section
\ref{sec:methodology}, after introducing notation and giving necessary
definitions, we outline a general methodology for constructing a
discrete finite element method that preserves the Lie point symmetries
of an ODE of arbitrary order. In Section \ref{sec:examples} we
examine several examples highlighting the intricacies of the
methodology.  In Section \ref{sec:numerics} we provide
numerical experiments that highlight favourable aspects of
symmetry-preserving finite element schemes.  In particular, our
simulations show that for certain equations the invariant schemes are
more accurate long term and can be implemented on larger elements than
their non-invariant counterparts.  Finally, in Section
\ref{sec:comparison} we contrast the ideas introduced in this paper
with those first considered in \cite{BihloValiquette:2018}.  In
particular we notice that our methodology extends the work of
\cite{BihloValiquette:2018} to non-projectable symmetry actions and is
much easier to generalise to higher order ODEs.

\revise{Before diving into the subject, we note that the present paper should be regarded as an exploratory study in assessing the feasibility of systematically incorporating symmetries into the finite element method.  Since this paper shows that this is possible, the next natural step will be to consider applications to partial differential equations, which offers a much richer playing ground.}

%%%%%
\section{General methodology} \label{sec:methodology}
%%%%%

Let 
\begin{equation} \label{eqn:ode}
  \ode{t,y,y_t,y_{tt},...,y_{t^{m+1}}} = 0,
\end{equation}
be an $\bc{m{+}1}$-th order ODE with initial conditions
\begin{equation}\label{eqn:initial_condition}
y(0)=y_0,\qquad \ldots,\qquad y_{t^m}(0)=y_{m}.
\end{equation}
To simplify the notation, we introduce the $\bc{m{+}1}$-th jet notation
\begin{equation}\label{eqn:jet}
y^{\bc{m+1}}=\bc{y,y_t,...,y_{t^{m+1}}},
\end{equation}
which collects the dependent variable
$y$ and its time derivatives up to order $m+1$.  The tuple $(t,y^{\bc{m+1}})$ provides local coordinates for the $\bc{m{+}1}$-th order
extended jet bundle $\jet{m+1}{}(\mathbb{R}^2,1)$ of
curves in the plane, \cite{Olver:1993}.
 
We assume that equation \eqref{eqn:ode} admits a certain Lie group of
point symmetries, which are \revise{transformations depending on the independent and dependent variables that map solutions of
\eqref{eqn:ode} to themselves}. There are many excellent textbooks on this
topic, such as \cite{BlumanAnco:2002,Olver:1993,Hydon:2000}.

\begin{definition}
  A Lie group $G$ is said to be a \emph{symmetry group} of equation
  \eqref{eqn:ode} if and only if for all $g\in G$ near the identity,
  the equality
  \begin{equation}
    \restrict{\ode{g\cdot t, g \cdot y^{\bc{m+1}}}}{\ode{%
        t,y^{\bc{m+1}}}=0}
    =0
  \end{equation}
  holds. 
\end{definition}

The symmetry group (or at least the infinitesimal symmetry generators)
of a given differential equation can easily be computed using symbolic
software packages such as {\sc Maple, Mathematica}, or {\sc Sage}.
Given the action of $G$ on the independent and dependent variables
\[
\g{t} : = g\cdot t,\qquad \g{y} : = g\cdot y, 
\]
the prolonged action on the time derivatives is obtained by implicit differentiation
\begin{equation}\label{eqn:prolonged action}
\g{y_{t^i}} = g\cdot y_{t^i} = \frac{\d^{i}\g{y}}{\d \g{t}^{\,i}}.
\end{equation}

In the sequel, we rewrite the $\bc{m{+}1}$-th order ODE \eqref{eqn:ode}
as a system of $m{+}1$ first order differential equations by
introducing the variables
\begin{equation} \label{eqn:aux}
  \vec{u} = (u_0,u_1,u_2,\ldots,u_m) = (y,y_t,\ldots,y_{t^m}),
\end{equation}
to obtain the system
\begin{equation} \label{eqn:sys}
  \begin{split}
    \ode{t,\vec{u},{u_m}_t} & = 0, \\
    {u_0}_t - u_1 & = 0, \\
    \vdots & \\
    {u_{m-1}}_t - u_{m} & = 0,
  \end{split}
\end{equation}
subject to the initial data 
\[
\vec{u}(0) = (y_0,y_1,\ldots,y_m).
\]
In light of how $\vec{u}$ is defined in \eqref{eqn:aux}, the induced
action of $G$ on $\vec{u}$ is given by the prolonged action
\eqref{eqn:prolonged action} after a change of variables.

\begin{remark}
  
  In general, the system of equations \eqref{eqn:sys} may possess more
  symmetries than the original equation \eqref{eqn:ode}. In the
  following, when working with system \eqref{eqn:sys}, we restrict
  ourselves to the symmetry group \revise{of} the original
  equation \eqref{eqn:ode}.
  
\end{remark}

\begin{example}[{\bf A \revise{working example}}]\label{ex:painleve1}

  To illustrate the notions introduced throughout this section, we
  will consider the nonlinear \revise{differential equation}
  \begin{equation} \label{eqn:painleve}
    y_{tt} - y^{-1} y_t^2 = 0,
  \end{equation}
  subject to the initial conditions
  \begin{equation}
    y(0)=y_0, \qquad y_t(0) = y_1,
  \end{equation}
  for some prescribed constants $y_0, y_1$, and with the assumption
  that $y(t) \neq 0$ for all $t \in \bs{0,T}$, for some $T>0$.  The
  differential equation \eqref{eqn:painleve} admits a six-parameter
  symmetry group of projectable transformations given by
  \[
    \g{t}=\frac{\alpha t + \beta}{\gamma t+\delta},\qquad \g{y} = (y^\lambda e^{ax+b})^{1/(\gamma x+\delta)},
  \]
  where $\alpha\delta-\beta\gamma = 1$, $a,b\in \mathbb{R}$, and
  $\lambda>0$.  In the following we will restrict ourselves to the
  two-parameter subgroup
  \begin{equation} \label{eqn:painleve-sym}
    \g{t} = t,\qquad \g{y} = \exp{at+b}\, y.
  \end{equation}
  This makes the exposition easier to follow by keeping the formulas
  simple, and furthermore we will see in Section \ref{sec:numerics}
  that the finite element method that preserves this two-parameter
  symmetry group yields desirable numerical results.  Up to order two, the
  prolonged action of \eqref{eqn:painleve-sym} is simply obtained by
  computing the product rule for differentiation:
  \begin{equation}
    \begin{aligned}
      \g{y_t} &= \frac{\d\g{y}}{\d t} = \frac{\d }{\d t}[\exp{at+b}y]=(ay+y_t)\exp{at+b},\\
      \g{y_{tt}} &= \frac{\d\g{y_t}}{\d t} = \frac{\d}{\d t}[(ay+y_t)\exp{at+b}] = (y_{tt}+2ay_t+a^2y)\exp{at+b}.
    \end{aligned}
  \end{equation}
  Then 
  \[
    \g{y_{tt}}-(\g{y})^{-1}(\g{y_t})^2=\big(y_{tt} - y^{-1}y_t^2\big)\exp{at+b} = 0,
  \]
  provided that $y$ is a solution of \eqref{eqn:painleve}, which
  confirms that \eqref{eqn:painleve-sym} is a symmetry (sub-)group of
  the equation.

Equation \eqref{eqn:painleve} is written as a system of first order ODE by introducing the variables
\begin{equation}\label{eqn:u-var}
\u = y,\qquad \v=y_t,
\end{equation}
so that \eqref{eqn:painleve} becomes
\begin{equation} \label{eqn:painleve2}
  \begin{split}
    \v_t - \u^{-1}\v^2 & = 0, \\
    \u_t - \v & = 0, \\
    \u(0) = y_0, \quad  \v&(0) = y_1.
  \end{split}
\end{equation}
Making the change of variables \eqref{eqn:u-var} into the transformation formulas \eqref{eqn:painleve-sym}, we conclude that the system of equations \eqref{eqn:painleve2} is invariant under the transformation group
\begin{equation} \label{eq:painleve-u-action}
\begin{aligned}
\g{t} &= t, &\qquad 
\g{\u} &= \exp{at+b}\u,\\
\g{\v} &= (a\u+\v)\exp{at+b},&\qquad
\g{\v_t} &= (\v_t+2a\v+a^2\u)\exp{at+b}.
\end{aligned}
\end{equation}
\end{example}

\subsection{Temporal finite element methods}
%%% Finite element stuff %%%
Let $T>0$, and consider a partition of the temporal interval
$\bs{0,T}$ into $N$ sub-intervals given by
$0=t_0 < t_1 < ... < t_N = T$.  For $n \in \{0,\ldots,N-1\}$, let
$I_n = ( t_n,t_{n+1} ]$ denote the \emph{finite element} of length
$\tau_n = t_{n+1}-t_n$.  Given a continuous or differentiable function
$y(t)$, we systematically use capital letters to denote the functions
finite element approximation. For example, $Y=Y(t)$ represents the
finite element function approximating $y(t)$. When there is no
ambiguity, we shall not explicitly write the time dependency of
functions.

\begin{definition}\label{def:tfes}

  Let $\bpoly{q}(I_n)$ denote the space of polynomials of degree $q$ on
  the element $I_n \subset [0,T]$.  Then the \emph{discontinuous finite
  element space} is
  \begin{equation} \label{def:dpoly}
    \dpoly{q} \bc{[0,T]} 
    =
    \{ Y : \left. Y \right|_{I_n} 
    \in \bpoly{q} \bc{I_n}, n=0,...,N-1 \}.
  \end{equation}
  The \emph{continuous finite element space} is defined
  analogously with global continuity enforced, i.e.
  \begin{equation} \label{def:cpoly}
    \cpoly{q} \bc{\bs{0,T}}
    =
    \dpoly{q} \bc{\bs{0,T}} \cap \mathcal{C}^0\bc{\bs{0,T}}
    ,
  \end{equation}
  where $\mathcal{C}^0\bc{\bs{0,T}}$ is the space of continuous
  functions over the interval $[0,T]$.
  
\end{definition}

We let $\cpoly{q} \bc{I_n}$ represent the localisation
of the continuous finite element space to a single element. Here, the
values of functions at time $t_n$ in $\cpoly{q}\bc{I_n}$ are fixed by
the right endpoint of the corresponding functions in
$\cpoly{q}\bc{I_{n-1}}$, for $\cpoly{q}\bc{I_0}$ the initial function
values are enforced by the initial conditions \eqref{eqn:initial_condition}. We note that
$\dpoly{q}\bc{I_n} \subset\mathcal{C}^{q}\bc{I_n}$ and
$\cpoly{q}\bc{I_ n} \subset\mathcal{C}^{q}\bc{I_n}$, i.e., both
discontinuous and continuous finite element approximations are smooth
over a single element. This property will be crucial in the sequel.

As we are restricting ourselves to initial value problems, we wish for
our finite element approximations to boast a time stepping
implementation. Moreover, as we have introduced the auxiliary
variables \eqref{eqn:aux}, we rewrite the ODE \eqref{eqn:ode} as the
first order system \eqref{eqn:sys} with dependent variable
$\vec{u}$. This allows for the design of finite element approximations
which encompass all initial value problems for ODEs.  As mentioned in
the introduction, we consider the \emph{continuous Galerkin method},
which is well studied for non-degenerate temporal ODEs,
\cite{EstepFrench:1994}. For such problems the method is described by
seeking for $\vec{U} \in \bc{\cpoly{q+1}\bc{I_n}}^{m+1}$, where
$n=\{0,\ldots,N-1\}$, such that
\begin{equation} \label{eqn:cg}
  \begin{split}
    \int_{I_n} \ode{t,\vec{U},{U_m}_t} \phi \di{t}
    & = 0
    \qquad \forall \phi \in \dpoly{q}\bc{I_n},
    \\
    \int_{I_n} \bc{{U_0}_t - U_1} \psi \di{t}
    & = 0
    \qquad \forall \psi \in \dpoly{q}\bc{I_n},
    \\
    \vdots &
    \\
    \int_{I_n} \bc{ {U_{m-1}}_t - U_m } \chi \di{t}
    & = 0
    \qquad \forall \chi \in \dpoly{q}\bc{I_n}
    ,
  \end{split}
\end{equation}
where $\vec{U}\bc{0}$ is fixed by the initial condition
$\vec{u}\bc{0}$, and $\vec{U}\bc{t_n}$ is determined by the solution
at the right endpoint on the previous element $I_{n-1}$. For a
detailed discussion and analysis of this method see
\cite{EstepFrench:1994, FrenchSchaeffer:1990}. The functions $\vec{U}$
are typically referred to as trial functions, and $\phi,\psi,..,\chi$
are known as test functions. We note that the upwind discontinuous
Galerkin approximation, \cite{Estep:1995}, is another widely used
finite element method for time evolving ODEs.  However, this method
does not conserve any underlying structures associated to the
continuous problem, \cite{EstepStuart:2002}. In this work, we will
refer to the continuous Galerkin formulation as ``standard''
\revise{in view of the following remark.}

\revise{
  \begin{remark}[{\bf Preservation of geometric structures by
    \eqref{eqn:cg}}] \label{rem:geometric}

    The standard continuous Galerkin formulation \eqref{eqn:cg} is
    well studied and understood. In fact, for a large class of
    geometrically interesting ODEs, namely Hamiltonian ODEs, it
    \emph{exactly} preserves the associated energy at the points
    $t_n$, see \cite{Hansbo:2001, self:thesis}. For an introduction to
    Hamiltonian ODEs and their geometric properties see
    \cite{HairerLubichWanner:2006, LeimkuhlerReich:2004}. In the
    lowest order case, i.e., when $q=0$, the standard finite element
    method is equivalent to the average vector field (AVF) discrete
    gradient method \cite{McLachlanQuispelRobidoux:1999}, and may therefore be
    viewed as a generalisation of this method. Furthermore, we note
    the close relation between the standard method and well known
    \emph{symplectic} integrators. Through an appropriate choice of
    Gaussian quadrature, i.e., choosing the number of quadrature points
    equal to the polynomial degree of the numerical solution, we
    obtain the $q+1$ point Gauss collocation method. This collocation
    method is, in turn, equivalent to a member of the Gauss family of
    Runge-Kutta methods which are known to be symplectic. For more
    information see \cite[\S 3.1]{self:thesis}. With this in mind, we
    expect the ``standard'' formulation to perform well numerically.

    We note that throughout this work we shall not study this method
    under quadrature approximation, instead evaluating integrals
    exactly where possible, or employing a quadrature which has a
    negligible error. This choice will be fully justified in
    Remark \ref{rem:quadcost}.
      
  \end{remark}
}
As the next example shows, the approximate weak formulation \eqref{eqn:cg} will, in general, not preserve the symmetries of the original system of differential equations \eqref{eqn:aux}.  

\begin{example}\label{ex:painleve2}
Continuing Example \ref{ex:painleve1}, the ``standard'' finite element formulation of the first order system \eqref{eqn:painleve2} is obtained by seeking for 
$\U,\V \in \cpoly{q+1}\bc{I_n}$ such that
\begin{equation} \label{eqn:painlevefe}
  \begin{split}
    \int_{I_n} \bc{
      \V_t - \U^{-1} \V^2
    } \phi \di{t}
    & = 0
    \qquad \forall \phi \in \dpoly{q}\bc{I_n},
    \\
    \int_{I_n} \bc{
      \U_t - \V
    } \psi \di{t}
    & = 0
    \qquad \forall \psi \in \dpoly{q}\bc{I_n}
    .
  \end{split}
\end{equation}
Observe that if $\U(t) = 0$ for some $t \in I_n$, then this finite
element is not well defined, so we must assume that this is not the
case numerically. However, we do not have an analytic assurance of
this requirement. Due to the local smoothness of the approximation,
the group action on $(t,U_0,U_1,{U_1}_t)$ is obtain by simply
replacing the lower case $u$ in \eqref{eq:painleve-u-action} by its
capital case counterpart $U$. Thus,
\begin{equation}\label{eqn:painleve2:symU}
\begin{aligned}
\g{t} &= t, &\qquad
\g{\U} &= \exp{at+b} \U,\\
  \g{\V} &= (a\U + \V)\exp{at+b},&\qquad
  \g{\V_t} &= (\V_t+2a\V+a^2\U)\exp{at+b}.
  \end{aligned}
\end{equation}
Therefore, the transformed weak form is 
\begin{equation}\label{eqn:transformed:painlevefe}
  \begin{split}
    \int_{I_n} \exp{at+b} \bc{
      \V_t - \U^{-1}\V^2 + a \U_t - a \V
    } \phi \di{t}
    & = 0, \\
    \int_{I_n} \exp{at+b} \bc{
      \U_t - \V
    } \psi \di{t}
    & = 0
    ,
  \end{split}
\end{equation}
observing that $\phi$ and $\psi$ are both arbitrary functions in the
same function space.  The transformed finite element method reduces to
\begin{equation} \label{eqn:painlevefe2}
  \begin{split}
    \int_{I_n} \exp{at+b} \bc{
      \V_t - \U^{-1} \V^2
    } \phi \di{t}
    & = 0, \\
    \int_{I_n} \exp{at+b} \bc{
      \U_t - \V
    } \psi \di{t}
    & = 0
    .
  \end{split}
\end{equation}
As the exponential function does not commute with the integral, we
conclude that the standard finite element method is not invariant
under the group action \eqref{eqn:painleve2:symU}. 
\end{example}

To construct finite element methods that preserve Lie point symmetries, we will use the method of equivariant moving frames, which we now summarise in the particular context of our problem.

\subsection{Moving frames and invariantisation} \label{sec:mfi}

For an introduction to the method of equivariant moving frames, we refer the reader to the original papers \cite{FelsOlver:1999,MariMansfield:2018,OlverPohjanpelto:2008,Olver:2019} and the textbook \cite{Mansfield:2010}.

Let $G$ be the symmetry group of the differential equation \eqref{eqn:ode}, or possibly a symmetry subgroup. Via the prolonged action, it induces an action 
\begin{equation}\label{eqn:u action}
\widehat{\vec{u}} = g\cdot \vec{u}
\end{equation}
on the variables \eqref{eqn:aux}, and remains a symmetry group of \eqref{eqn:sys}.  Since the function $\vec{U}=(U_0,\ldots,U_m)$ is smooth in the interior of any element $I_n$, the action \eqref{eqn:u action} also applies to $\vec{U}$ inside $I_n$, to give $\widehat{\vec{U}} = g\cdot \vec{U}$.  It is prolonged to 
\[
\vec{U}^{\bc{1}} = (U_0,\ldots,U_m,{U_0}_t,\ldots, {U_m}_t)
\]
using the usual prolongation formula
\[
\widehat{{U_i}_t} = \frac{d\widehat{U}_i}{d\widehat{t}},\qquad i=0,\ldots,m.
\]
Therefore, in the following, we consider the action of the Lie group $G$ on the first order jet space $\jet{1}{I_n}=\mathrm{J}^{(1)}(\mathbb{R}^{m+1},1)\big|_{I_n}$ with local coordinates \revise{$\vec{Z}^{(1)}  = (t,\vec{U}^{(1)})$. We combine the Lie group and the first order jet space by introducing the \emph{first order lifted bundle} $\mathcal{B}^{(1)}_{I_n} = G\times \jet{1}{I_n}$.  This bundle admits a Lie groupoid structure, \cite{Mackenzie:2005}, with source map $\source(g,\vec{Z}^{(1)}) = \vec{Z}^{(1)}$ corresponding to the projection onto $\jet{1}{I_n}$ and the target map $\target(g,\vec{Z}^{(1)}) = \widehat{\vec{Z}}{}^{(1)}=g\cdot\vec{Z}^{(1)}$ given by the prolonged group action.}

\begin{definition}
  
  Let $H\colon \jet{1}{I_n} \to \mathbb{R}$ be a function.  The
  \emph{lift} of $H$ is the function
  \[
    \bl(H(\vec{Z}^{(1)})) = \revise{\target^*H(\vec{Z}^{(1)})=H(\widehat{\vec{Z}}{}^{\bc{1}})}
  \]
  obtained by \revise{substituting the arguments $\vec{Z}^{(1)}$ of $H$ by their prolonged action expressions $\widehat{\vec{Z}}{}^{\bc{1}}$.}
\end{definition}

The lift of a function is a new function defined on the \revise{lifted
  bundle $\mB^{(1)}_{I_n}$}. The lift map
$\bl$ extends to differential forms, see \cite{KoganOlver:2003} for more detail. The (horizontal) lift of
$\di{t}$ is the one-form
\begin{equation} \label{eqn:omega}
  \omega =  \bl(\mathrm{d}t) = \bigg(\frac{\partial \g{t}}{\partial t} + \frac{\partial \g{t}}{\partial U_0}\cdot U_1\bigg)\di{t},
\end{equation}
obtained by taking the (horizontal) differential of $\widehat{t} = g\cdot t$, where the group parameters $g=(g^1,\ldots, g^r)$ are treated as constants.  In order to have a well defined action of the symmetry group $G$ on the system of equations \eqref{eqn:cg}, there remains to define the action of $G$ on the basis functions $\phi \in \dpoly{q}\bc{I_n}$.  By definition, the function $\phi$ is a degree $q$ Lagrangian interpolating function depending on the nodes of the element $I_n$, i.e., $t_n$, $t_{n+1}$, a finite number of interior points $t_i\in (t_n,t_{n+1})$, and obviously the continuous time variable $t$.   We write these dependencies explicitly as
\[
\phi = \phi(t_n,\ldots,t_{n+1};t),
\]
where in $t_n,\ldots,t_{n+1}$, we collect the dependency of
$\phi$ on the nodes $t_n$, $t_{n+1}$ and the interior points
$t_i$. Then,
the lift of $\phi$ is defined as the function
\begin{equation}\label{eqn:lift-phi}
\bl(\phi) = \widehat{\phi} = \revise{\target^* \phi} =  \phi(g\cdot t_n,\ldots,g \cdot t_{n+1};g\cdot t)
\end{equation}
where the group $G$ acts via the product action on $t_n$, $\ldots$, $t_{n+1}$, and $t$, simultaneously.

To simplify the exposition, we now introduce the notation
\begin{equation} \label{eqn:function}
  \fe{\revise{\vec{Z}^{(1)}},\bphi,\mathrm{d}t} = 0,
\end{equation}
to refer to the system of equations \eqref{eqn:cg}, where $\bphi=(\phi,\psi,\ldots,\chi)$ denotes the test functions occurring in \eqref{eqn:cg}.  Then the lift of \eqref{eqn:cg}, in other words, the action of $G$ on \eqref{eqn:cg} is defined as 
\begin{align*}
\fe{\bl(\revise{\vec{Z}^{(1)}},\bphi,\mathrm{d}t)} &= \fe{\revise{\bl(\vec{Z}^{(1)})},\bl(\bphi),\bl(\mathrm{d}t)}\\
&=\fe{\revise{g\cdot \vec{Z}^{(1)}},\bphi(g\cdot t_n,\ldots,g\cdot t_{n+1};g\cdot t),\omega}
                                                                                                              ,
\end{align*}
where $\omega$ is given in \eqref{eqn:omega}.
\begin{definition}
The system of equations \eqref{eqn:cg} is said to be \emph{invariant} under the action of a Lie group $G$ if and only if
\begin{equation}
  \bl(\fe{\revise{\vec{Z}^{(1)}},\bphi,\mathrm{d}t})\big|_{\fe{\revise{\vec{Z}^{(1)}},\bphi,\mathrm{d}t}=0}=0.
\end{equation}
\end{definition}

As mentioned in the previous section, the approximate weak form \eqref{eqn:cg} will not, in general, preserve the symmetries of the original system \eqref{eqn:aux}.  To obtain a weak formulation that preserves the symmetries, we follow the methodology of invariantisation, which is based on the introduction of a moving frame. 

\begin{definition} \label{def:movingframe}

Let $G$ be a Lie group acting on $\jet{1}{I_n}$.  A \emph{right moving
  frame} is a $G$-equivariant map $\rho\colon \jet{1}{I_n} \to G$ satisfying the equality
  \begin{equation} \label{eqn:movingframe}
    \mf{g\cdot \vec{Z}^{\bc{1}}}
    =
    \mf{\vec{Z}^{\bc{1}}} \cdot g^{-1}
    \qquad \text{with}\qquad g\in G.
  \end{equation}
    
\end{definition}

To guarantee the existence of a moving frame, we need to impose two regularity assumptions on the group action, \cite{FelsOlver:1999}.

\begin{definition}

  The Lie group $G$ is said to act \emph{freely} at $\vec{Z}^{\bc{1}}\in \jet{1}{I_n}$ if the isotropy group
  \begin{equation}
    G_{\vec{Z}^{\bc{1}}}
    =
    \big\{ g \in G \textrm{ } | \textrm{ } g \cdot \vec{Z}^{\bc{1}} = \vec{Z}^{\bc{1}} \big\} = \{e\}
  \end{equation}
  is trivial. The group action is said to be locally free at $\vec{Z}^{\bc{1}}$ if $G_{\vec{Z}^{\bc{1}}}$ is
  a discrete subgroup.  We say that $G$ acts (locally) freely on $\jet{1}{I_n}$ if the action of $G$ is (locally) free at every point
  $\vec{Z}^{\bc{1}} \in \jet{1}{I_n}$.

\end{definition}
  
\begin{definition}
  
  The action of $G$ on $\jet{1}{I_n}$ is said to be \emph{regular} if the
  orbits form a regular foliation, that is to say that the orbits are
  submanifolds of dimension $r=\dim G$.%, see \cite[Definition
  %1.25]{Olver:1995}.

\end{definition}

\begin{theorem}

  Let $G$ act (locally) freely and regularly on $\jet{1}{I_n}$, then a moving frame
 exists in a neighbourhood of every point $\vec{Z}^{\bc{1}} \in \jet{1}{I_n}$.
  
\end{theorem}

The construction of a moving frame is based on the introduction of a cross-section.

\begin{definition}

A \emph{cross-section} $\mathcal{K} \subset \jet{1}{I_n}$ is a submanifold of co-dimension $r=\dim G$ transverse to the group orbits.

\end{definition}

Given a cross-section $\mathcal{K} \subset \jet{1}{I_n}$, the right moving frame at the point $\vec{Z}^{(1)}=\big(t,\vec{U}^{\bc{1}}\big)$ is defined as 
unique group element $\rho\in G$ such that
\begin{equation}\label{eqn:mfK}
\mf{\vec{Z}^{\bc{1}}}\cdot \vec{Z}^{\bc{1}} \in \mathcal{K}.
\end{equation}
In practice, a moving frame is frequently obtained by choosing a
coordinate cross-section, which consists of setting certain
coordinates of the jet $\vec{Z}^{\bc{1}}=(t,\vec{U}^{(1)})$
equal to constant values.  Due to the definition of the variables
$\vec{u}$ in \eqref{eqn:aux}, and the fact that the transformation
group $G$ comes from the symmetry group of the differential equation
\eqref{eqn:ode}, we work under the assumption that it is possible to
determine a coordinate cross-section by setting $r=\dim G$ coordinates
from $\vec{Z}^{(0)}=(Z_0,\ldots,Z_{m+1}) = (t,\vec{U}^{(0)})$ to
constant values.  Thus, let
\begin{equation} \label{eqn:cs}
  \cs = \big\{Z_{i_1} = c_{i_1},\ldots, Z_{i_r}=c_{i_r}\big\}\subset \jet{1}{I_n},
 \end{equation} 
 where $0\leq i_1<i_2<\cdots <i_r\leq m+1$.  The requirement \eqref{eqn:mfK} leads to the \emph{normalisation equations}
\begin{equation} \label{eqn:normal}
  g \cdot Z_{i_1} = c_{i_1},\qquad \ldots,\qquad g\cdot Z_{i_r} = c_{i_r}.
\end{equation}
Solving for the group parameters $g$ yields the moving frame
$g=\rho(\vec{Z}^{(1)})=\rho(\vec{Z}^{(0)})$.

\begin{example}\label{ex:painleve3}
Continuing Example \ref{ex:painleve2}, we now construct a moving frame for the action \eqref{eqn:painleve2:symU}.  Working on the set where $U_0\neq 0$, the action \eqref{eqn:painleve2:symU} is free and regular, therefore a moving frame exists.  A coordinate cross-section  to the group orbits is given by
\begin{equation} \label{eqn:painleve:crosssection}
  \cs = \bw{\U=\mathrm{sign}\bc{\U},\;\V=0}.
\end{equation}
The corresponding normalisation equations are
\begin{equation}
  \exp{at+b} \U = \mathrm{sign}\bc{\U}
  ,
  \qquad
  (a\U+\V)\exp{at+b} = 0
  .
\end{equation}
Solving these two equations for the group parameters $a,b$ leads to the right moving frame
\begin{equation}\label{eqn:painlevemf}
  a = - \frac{\V}{\U}
  ,
  \qquad
  b = t \cdot \frac{\V}{\U} - \ln{\abs{\U}}
  .
\end{equation}
\end{example}

Given a moving frame, there is a systematic procedure for mapping non-invariant quantities (e.g.\ functions, differential forms, functionals) to their invariant analogue, known as \emph{invariantisation}.  To introduce the invariantisation map, we note that a moving frame $\rho\colon \jet{1}{I_n} \to G$ induces a moving frame section $\varrho\colon \jet{1}{I_n} \to \mB^{(1)}_{I_n}$ in the lifted bundle $\mB^{(1)}_{I_n}$ defined by
\[
\varrho(\vec{Z}^{(1)}) = (\vec{Z}^{(1)},\rho(\vec{Z}^{(1)})).
\]
Then, given a function $H\colon \mB^{(1)}_{I_n}\to \mathbb{R}$, or more generally a differential form, we can consider the moving frame pull-back
\[
\varrho^*[H(\vec{Z}^{\bc{1}},g)]=H(\vec{Z}^{\bc{1}},\rho(\vec{Z}^{(1)})),
\]
obtained by replacing the group element $g$ by the moving frame $\rho(\vec{Z}^{(1)})$.

\begin{definition}
Let $\rho\colon \jet{1}{I_n} \to G$ be a moving frame.  The \emph{invariantisation} of the system of equations \eqref{eqn:cg} is defined as
\begin{equation}\label{eqn:invariantisation}
0 =\iota[\fe{\revise{\vec{Z}^{(1)}},\bphi,\mathrm{d}t}] = \fe{\varrho^* [\bl(\revise{\vec{Z}^{(1)}},\bphi,\mathrm{d}t)]},
\end{equation}
where we define $\varrho^*[\bl(\bphi)]=\bphi(\rho\cdot t_n,\ldots,\rho\cdot t_{n+1};\rho\cdot t)$.
\end{definition}

Therefore, the system \eqref{eqn:invariantisation} is obtained by first lifting its arguments, i.e.\ by acting on the arguments with the symmetry group $G$, and then replacing the group element $g$ by the moving frame $\rho$.  The map $\iota$ is called the \emph{invariantisation map}.  The fact that \eqref{eqn:invariantisation} is invariant follows from the $G$-equivariance of the moving frame $\rho$.  We note that if the finite element approximation
is already invariant, then the invariantisation of \eqref{eqn:invariantisation} will simply yield back the original system of equations.

\begin{example}\label{ex:painleve4}
As observed in Example \ref{ex:painleve2}, the discrete weak
formulation \eqref{eqn:painlevefe} is not invariant under the group of
transformations \eqref{eqn:painleve2:symU}.  To obtain a weak
formulation that will remain invariant under the group action, it
suffices to invariantise \eqref{eqn:painlevefe}.  The lift of the
discrete weak form has already been computed in
\eqref{eqn:painlevefe2}.  Using the moving frame
\eqref{eqn:painlevemf} computed in Example \ref{ex:painleve3}, the
moving frame pull-back of \eqref{eqn:painlevefe2} is simply obtained
by substituting the group normalisations \eqref{eqn:painlevemf} into
\eqref{eqn:painlevefe2}.  The result is the symmetry-preserving finite element scheme
\begin{equation} \label{eqn:painleve:invariant}
  \begin{split}
      \int_{I_n} \U^{-1}\bc{
      \V_t - \U^{-1} \V^2
    } \phi \di{t}
    & = 0
    \qquad \forall \phi \in \dpoly{q}\bc{I_n},
    \\
    \int_{I_n} \U^{-1}\bc{
      \U_t - \V
    } \psi \di{t}
    & = 0
    \qquad \forall \psi \in \dpoly{q}\bc{I_n}
    .
  \end{split}
\end{equation}
The invariance of the resulting scheme may be verified simply by
applying the group action \eqref{eqn:painleve2:symU}.
\end{example}

To summarise the previous two sections, we now carefully discuss the invariantisation procedure of a finite element approximation and its various intricacies.  More examples of the 
invariantisation procedure are considered in Section \ref{sec:examples}.

\subsection{Outline of the methodology for constructing symmetry-preserving finite element methods} \label{sec:genmeth}

To simplify the presentation, we restrict our exposition to third order ODEs, which encapsulates all examples considered in Section \ref{sec:examples}.  Following the exposition introduced in the previous sections, it is straightforward to extend the methodology to higher order equations.  

Our starting point is a third order ODE
\begin{equation}\label{eqn:3ODE}
\ode{t,y,y_t,y_{tt},y_{ttt}}=0,
\end{equation}
with initial condition $y(0)=y_0$, $y_t(0)=y_1$, $y_{tt}(0) = y_2$, which admits a symmetry group $G$.  The 
 steps for constructing a symmetry-preserving continuous Galerkin finite element scheme are as follows:

\subsubsection*{Step 1}
Introduce the auxiliary variables
\begin{equation} \label{eqn:g:aux}
  u_0=y,
  \qquad
  u_1 = {u_0}_t=y_t,
  \qquad
  u_2 = {u_1}_t=y_{tt},
\end{equation}
and recast the third order equation \eqref{eqn:3ODE} as a system of first order ODEs
\begin{equation} \label{eqn:g:sys}
  \begin{split}
    \ode{t,u_0,u_1,u_2,{u_2}_t} & = 0, \\
    {u_0}_t - u_1 & = 0, \\
    {u_1}_t - u_2 & = 0
    ,
  \end{split}
\end{equation}
with initial conditions $u_0(0)=y_0$, $u_1(0)=y_1$, $u_2(0)=y_2$.  We note that while
$\ode{t,u_0,u_1,u_2,{u_2}_t}$ is uniquely prescribed here, there may be other,
equally valid, ways to describe the system of equations \eqref{eqn:g:sys}.  Such an example is provided in Example \ref{example:b}.

From the definition of the variable $u$ in \eqref{eqn:g:aux}, the symmetry group $G$ induces the transformations
\[
\g{t} = g \cdot t,\qquad \g{u}_0 = g\cdot u_0,\qquad \g{u_1}=g\cdot u_1,\qquad 
\g{u_2} = g\cdot u_2
\]
via the prolonged action.  Furthermore,
\[
\g{{u_0}_t} = g\cdot {u_0}_t = \frac{d\g{u_0}}{d\g{t}},\qquad \g{{u_1}_t} = g\cdot {u_1}_t = \frac{d\g{u_1}}{d\g{t}}.
\]

\subsubsection*{Step 2}
Formulate a ``standard'' finite element approximation in the spirit of
\eqref{eqn:cg} by seeking $U_0,U_1,U_2 \in \cpoly{q+1}\bc{I_n}$ such that
\begin{equation} \label{eqn:g:fe}
  \begin{split}
    \int_{I_n} \ode{t,U_0,U_1,U_2,{U_2}_t} \phi \di{t} & = 0
    \qquad
    \forall \phi \in \dpoly{q}\bc{I_n},
    \\
    \int_{I_n} \bc{{U_0}_t - U_1} \psi \di{t} & = 0
    \qquad
    \forall \psi \in \dpoly{q}\bc{I_n},
    \\
    \int_{I_n} \bc{{U_1}_t - U_2} \chi \di{t} & = 0
    \qquad
    \forall \chi \in \dpoly{q}\bc{I_n}
    ,
  \end{split}
\end{equation}
where $U_0(t_n),U_1(t_n),U_2(t_n)$ are fixed by either the solution on the
previous element or the initial data.

\subsubsection*{Step 2.1}

This step is not necessary, but if possible we wish for our space of
test functions $\dpoly{q}\bc{I_n}$ to be invariant under the group
action $\g{t} = g\cdot t$. If this is not the case, it may be possible
to choose a more appropriate function space for a given symmetry group
action, however, no such examples are presented in the sequel.

\subsubsection*{Step 3}

Construct a moving frame.  Assuming the action of $G$ on 
\begin{equation}
  \vec{Z}^{\bc{0}} := \bc{t,\U,\V,\W}
  ,
\end{equation}
is free and regular, choose a cross-section
\begin{equation} \label{eqn:g:cs}
  \cs = \bw{Z_{i_1} = c_{i_1},\ldots,Z_{i_r}=c_{i_r}},\qquad \text{where}\qquad r=\dim G \leq 4,
\end{equation}
and solve the normalisation equations $g\cdot Z_{i_1} = c_{i_1},\ldots, g\cdot Z_{i_r}=c_{i_r}$ for the group parameters to obtain the moving frame $\rho(\vec{Z}^{(0)})$.

\subsubsection*{Step 4}

Invariantise \eqref{eqn:g:fe} according to the invariantisation formula \eqref{eqn:invariantisation}, which is obtained by first lifting $(t,\vec{U}^{(1)})$, $\bphi=(\phi,\psi,\chi)$, and $\mathrm{d}t$ using the group action $G$ and then substituting the group elements my their moving frame expressions. The result is a discrete weak formulation of \eqref{eqn:g:sys} that is invariant under the action of the symmetry group of the differential equation \eqref{eqn:3ODE}.

\revise{
\begin{remark}[\bf Normalisation equations without exact solutions] \label{rem:unsolve}
Steps 3 and 4 are based on the assumption that the normalisation equations can be solved for the group parameters to obtain an explicit moving frame $\rho(\vec{Z}^{(0)})$, which is then used to invariantised \eqref{eqn:g:fe}.  If the normalisation equations cannot be solved analytically, it is still possible to obtain a symmetry-preserving finite element method as follows.  Instead of steps 3 and 4, consider the enlarged system of equations consisting of the normalisations with the lift of \eqref{eqn:cg}:
\begin{gather*}
g\cdot Z_{i_1} = c_{i_1},\ldots, g\cdot Z_{i_r}=c_{i_r},\\
\fe{g\cdot \vec{Z}^{(1)},\bphi(g\cdot t_n,\ldots,g\cdot t_{n+1};g\cdot t),\omega}=0.
\end{gather*}
Solving these equations numerically to a high precision will then constitute the desired invariant finite element method.
\end{remark}}

\revise{\begin{remark}[\bf Generalisation to systems of ODEs] \label{rem:sys}
In the above exposition, we have restricted our attention to a single ODE.  The procedure for constructing a symmetry-preserving finite element method can easily be extended to systems of ODEs.  From a notational standpoint, simply replace the original ODE \eqref{eqn:ode} by a system of $\ell$ ODEs and let $y(t)=(y^1(t),\ldots,y^\ell(t))$ denote a vector of dependent variables.  Then the rest of the exposition remains essentially unchanged, except that $u_i=y_{t^i}$ now denotes the time derivative of an $\ell$-dimensional vector and the discontinuous and continuous finite element spaces should be replaced by the Cartesian products
\begin{gather*}
\dpoly{q} \bc{[0,T]}^\ell
    =
    \{ (Y^1,\ldots,Y^\ell) : \left. Y^i \right|_{I_n} 
    \in \bpoly{q} \bc{I_n},\; i=1,\ldots,\ell,\; n=0,...,N-1 \},\\
(\cpoly{q} \bc{\bs{0,T}})^\ell = \dpoly{q} \bc{\bs{0,T}}^\ell \cap \mathcal{C}^0\bc{\bs{0,T}}^\ell.
\end{gather*}
\end{remark}}

\section{Examples} \label{sec:examples}

In this section we provide several examples highlighting different aspects of the methodology for constructing symmetry-preserving continuous Galerkin finite element methods.   \revise{Our first example is concerned with second order linear equations which occur in classical mechanics, electric circuits, and many other branches of science.  Our second example involves the Schwarzian differential equation which occurs in geometry, the theory of Sturm--Liouville equations, \cite{OvsienkoTabachnikov:2009}, and in the study of gravity--dilation--antisymmetric tensor systems, \cite{ZhouZhu:1999,ZhouZhu:1999-2}.  Finally, our last two examples have been selected so as to admit interesting symmetry groups and to illustrate certain aspects of our methodology.} In Section \ref{sec:numerics}, numerical simulations using the obtained schemes are performed and compared to standard non-invariant finite element methods.

\begin{example}[{\bf Second order linear ODE}] \label{example:glode}
Our first example is not an illustration of the invariantisation
procedure per say, but rather an illustration of the fact that for certain differential equations the continuous Galerkin method is naturally invariant.  To this end, let $y=y(t)$ satisfy the initial value problem
\begin{equation} \label{eqn:glode}
  \begin{gathered}
    y_{tt} + p(t) y_t + q(t) y = f(t), \\
    y(0)=y_0, \qquad y_t(0)=y_1
    ,
  \end{gathered}
\end{equation}
for some prescribed constants $y_0$ and $y_1$. This equation possesses a
two-parameter symmetry group given by
\begin{equation} \label{eqn:glode:sym}
  \g{t} = t,
  \qquad
  \g{y} = y + \epsilon_1 \alpha(t) + \epsilon_2 \gamma(t)
  ,
\end{equation}
where $\epsilon_1,\epsilon_2\in \mathbb{R}$, and  $\alpha(t)$, $\gamma(t)$ are linearly independent
solutions of the corresponding homogeneous equation
\begin{equation} \label{eqn:glode:homo}
  x_{tt} + p(t) x_t + q(t) x = 0
  .
\end{equation}
This symmetry group reflects the fact that the differential equation is linear. 
 Introducing the auxiliary variables $\u=y$ and $\v=y_t$,
 we may rewrite \eqref{eqn:glode} as the system of first order equations
\begin{equation} \label{eqn:glode2}
  \begin{split}
    \v_t + p(t) \v + q(t) \u & = f(t), \\
    \u_t - \v & = 0, \\
    \u(0) = y_0, \qquad  \v(0)& = y_1
    ,
  \end{split}
\end{equation}
which is invariant under the two-parameter group of transformations
\begin{equation} \label{eqn:glode2:sym}
  \begin{split}
    \G{t} = t,
    \qquad
    \G{\u} = \u + \epsilon_1 \alpha + \epsilon_2 \gamma,
    \qquad
    \G{\v}  = \v + \epsilon_1 \alpha_t + \epsilon_2 \gamma_t.
  \end{split}
\end{equation}
Applying the finite
element discretisation \eqref{eqn:cg} results in us seeking for
$\U,\V\in \cpoly{q+1}\bc{I_n}$ such that
\begin{equation} \label{eqn:glodefe}
  \begin{split}
    \int_{I_n} \bc{
      {\V}_t + p(t) \V + q(t) \U
      - f(t) }
    \phi \di{t}
    & = 0
    \qquad
    \forall \phi \in \dpoly{q}\bc{I_n},
    \\
    \int_{I_n} \bc{
      \U_t - \V }
    \psi \di{t}
    & = 0
    \qquad
    \forall \psi \in \dpoly{q}\bc{I_n}
    .
  \end{split}
\end{equation}
Since the independent variable $t$ is an invariant of the group action, it follows that $\di{t}$ and the test functions $\phi$ and $\psi$ remain unchanged under the group action \eqref{eqn:glode2:sym}.  Thus, the transformed finite element scheme is simply obtained by substituting the transformation rules
\[
    \G{\U} = \U + \epsilon_1 \alpha + \epsilon_2 \gamma,
    \qquad
    \G{\V}  = \V + \epsilon_1 \alpha_t + \epsilon_2 \gamma_t,
\]
and
\[
\g{\V_t} = \frac{d\g{\V}}{\di{t}} = \V_t + \epsilon_1 \alpha_{tt} + \epsilon_2 \gamma_{tt}
\]
into \eqref{eqn:glode2:sym}.  The result is
\begin{equation}
  \begin{gathered}
    \int_{I_n} \bc{
      \V_t + p(t) \V + q(t) \U -  f(t)
    } \phi
    + \epsilon_1 \bc{
      \alpha_{tt} + p(t) \alpha_t + q(t) \alpha
    } \phi
    + \epsilon_2 \bc{
      \gamma_{tt} + p(t) \gamma_t + q(t) \gamma
    } \phi
    \di{t}
     = 0,
    \\
    \int_{I_n} \bc{
     {U_0}_t - U_1
    } \psi \di{t}
     = 0
    .
  \end{gathered}
\end{equation}
Since $\alpha$ and $\gamma$ are solutions of the homogeneous
problem \eqref{eqn:glode:homo}, it follows that the terms in $\epsilon_1$ and $\epsilon_2$ vanish, and
therefore, the finite element scheme \eqref{eqn:glodefe} is invariant under the group action \eqref{eqn:glode2:sym} without any need to invariantise the functional.  In other words, since the scheme is linear it preserves the superposition principle.
\end{example}

\begin{example}[{\bf Schwarzian differential equation}] \label{example:order3}
Consider the third order ODE
\begin{equation} \label{eqn:order3}
  \frac{y_{ttt}}{y_t}
  - \frac{3}{2} \bigg(\frac{y_{tt}}{y_t}\bigg)^2
  = F(t)
  ,
\end{equation}
subject to the initial data $y(0)=y_0, y_t(0)=y_1,
y_{tt}(0)=y_2$, and where we assume that $y_t(t)\neq 0$ for all $t \in
[0,T]$. The differential equation is known, \cite{Mansfield:2010}, to be invariant under the linear 
fractional group action
\begin{equation} \label{eqn:order3:sym}
  \g{t} = t, \qquad
  \g{y} = \frac{\alpha y + \beta}{\gamma y + \delta}, \qquad
  \alpha \delta - \beta \gamma = 1
  .
\end{equation}
By letting $\u=y$,  we can rewrite \eqref{eqn:order3} as the system of first order differential equations
\begin{equation} \label{eqn:order3:sys}
  \begin{split}
    \frac{\w_t}{\v} - \frac{3}{2} \bigg(\frac{\w}{\v}\bigg)^2 & = F(t), \\
    \u_t - \v & = 0, \\
    \v_t - \w & = 0,
  \end{split}
\end{equation}
which is invariant under the extended group of transformations
\begin{equation} \label{eqn:order3:sym2}
  \begin{split}
    \g{t} = t, \qquad
    \g{\u} = \frac{\alpha \u + \beta}{\gamma \u + \delta}, \qquad
    \g{\v} = \frac{\v}{\bc{\gamma \u + \delta}^2}, \qquad
    \g{\w} = \frac{\w}{\bc{\gamma \u + \delta}^2}
    - \frac{2\gamma(\v)^2}{\bc{\gamma \u + \delta}^3}
    .
  \end{split}
\end{equation}
The ``standard'' finite element approximation is obtained
by seeking for $\U,\V,\W \in \cpoly{q+1}\bc{I_n}$ such
that
\begin{equation} \label{eqn:order3:fe}
  \begin{split}
    \int_{I_n} \Bigg(
      \frac{\W_t}{\V}
      - \frac32 \bigg(\frac{\W}{\V}\bigg)^2
      - F(t) \Bigg) \phi\di{t} & = 0
    \qquad \forall \phi \in \dpoly{q}\bc{I_n}, \\
    \int_{I_n} \bc{
      \U_t - \V } \psi \di{t} & = 0
    \qquad \forall \psi \in \dpoly{q}\bc{I_n}, \\
    \int_{I_n} \bc{
      \V_t - \W } \chi \di{t} & = 0
    \qquad \forall \chi \in \dpoly{q}\bc{I_n}
    ,
  \end{split}
\end{equation}
subject to appropriate initial data. Applying the group action
\eqref{eqn:order3:sym2} to \eqref{eqn:order3:fe}, we find the transformed
functionals
\begin{equation} \label{eqn:order3:fe2}
  \begin{gathered}
    \int_{I_n} \Bigg(
      \frac{\W_t}{\V}
      - \frac{3}{2} \bigg(\frac{\W}{\V}\bigg)^2
      - F(t) + \frac{ 6 \gamma^2 \bc{
        \U_t\V - \V^2}}{\bc{\gamma \U + \delta}^2} \
      + \frac{\gamma \bc{
        6 \V \W - 4 \V \V_t - 2 \U_t \W }}{\V\bc{\gamma \U + \delta}} 
    \Bigg)\phi \di{t}  = 0,
    \\
    \int_{I_n} 
      \Bigg(\frac{\U_t - \V}{\bc{\gamma \U + \delta}^2}\Bigg)
    \psi \di{t}  = 0,\qquad\qquad
    \int_{I_n} \Bigg(
      \frac{\V_t - \W}{\bc{\gamma \U + \delta}^2}
      + \frac{2\gamma \V^2 - 2\gamma \V \U_t}{
        \bc{\gamma \U + \delta}^3}\Bigg)
    \chi \di{t}  = 0
    ,
  \end{gathered}
\end{equation}
which shows that the weak formulation is not invariant.  We observe
that this finite element approximation is consistent, i.e.,
substituting sufficiently globally smooth trial and test functions
into the numerical approximation gives back the original system of
ODEs \eqref{eqn:order3:sys}.

A moving frame is obtained by choosing the cross-section
\begin{equation} \label{eqn:order3:cs}
  \cs = \bw{\U=0,\; \V=1,\;\W=0}.
\end{equation}
This yields the normalisation equations
\begin{equation} \label{eqn:order3:normal}
  \frac{\alpha \U + \beta}{\gamma \U + \delta} = 0, \qquad
  \frac{\V}{\bc{\gamma \U + \delta}^2} = 1, \qquad
  \frac{\W}{\bc{\gamma \U + \delta}^2}
  - \frac{2\gamma \V^2}{\bc{\gamma \U + \delta}^3} =0, \qquad
  \alpha \delta - \beta \gamma = 1
  ,
\end{equation}
the solution of which gives the moving frame
\begin{equation} \label{eqn:order3:mf}
  \alpha = \pm \V^{-\frac12}
  \qquad
  \beta = \mp \U \V^{-\frac12}
  \qquad
  \gamma = \pm \W \V^{-\frac32}
  \qquad
  \delta = \mp \U\W\V^{-\frac32} \pm \V^{\frac12}
  .
\end{equation}
The invariantisation of the weak form \eqref{eqn:order3:fe} is then obtained by substituting the
group normalisations \eqref{eqn:order3:mf} into the transformed functionals \eqref{eqn:order3:fe2}. The result is the symmetry-preserving finite element method
\begin{equation} \label{eqn:order3:fe:invariant}
  \begin{split}
    \int_{I_n} \Bigg(
      \frac{\W_t}{\V}
      - 2 \frac{\V_t\W}{\V^2}
      + \frac12 \frac{\U_t\W^2}{\V^3}
      - F(t)
    \Bigg) \phi \di{t} & = 0
    \qquad \forall \phi \in \dpoly{q}\bc{I_n}, \\
    \int_{I_n} 
      \Bigg(\frac{\U_t - \V}{\V}\Bigg)
     \psi \di{t} & = 0
    \qquad \forall \psi \in \dpoly{q}\bc{I_n}, \\
    \int_{I_n} \Bigg(
      \frac{\V_t - \W}{\V}
      + \frac{\W}{\V^3} \bc{ \V^2 - \V \U_t }
    \Bigg) \chi \di{t} & = 0
    \qquad \forall \chi \in \dpoly{q}\bc{I_n}
    .
  \end{split}
\end{equation}
\end{example}

\begin{example}[{\bf A second order quasi-linear ODE}] \label{example:b}
We now move our attention to the second order ODE
\begin{equation} \label{eqn:b}
  \begin{gathered}
    t^2y_{tt} + 4t y_t + 2 y = \bc{2ty + t^2y_t}^{\frac12}, \\
    y(1) = y_0, \quad y_t(1) = y_1
    ,
  \end{gathered}
\end{equation}
for some constants $y_0, y_1$.  This differential equation admits a two-parameter symmetry group 
with nontrivial temporal action given by
\begin{equation} \label{eqn:b:sym}
  \g{t} = \exp{a} t + b
  ,
  \qquad
  \g{y} = \frac{\exp{3a} t^2 y}{\bc{\exp{a} t + b}^2}
  ,
\end{equation}
where $a,b\in \mathbb{R}$.  When considering a finite domain in both the continuous and discrete
setting, it is important to note that due to dilation and translation in time, the action will accordingly
dilate and shift the domain of consideration. With this in mind, the
``invariant'' scheme we will consider for this problem will be
permitted to dilate and shift the elements through time. Introducing the auxiliary variable $\u=y$ we may
write \eqref{eqn:b} as
\begin{equation} \label{eqn:b2}
  \begin{split}
    t^2\v_{t} + 4t \u_t + 2 \u - \bc{2t\u + t^2\u_t}^{\frac12} & = 0,
    \\
    \u - \v_t & = 0
    ,
  \end{split}
\end{equation}
which is invariant under the extended group action
\begin{equation} \label{eqn:b2:sym}
  \g{t} = \exp{a} t + b
  ,
  \qquad
  \g{\u} = \frac{\exp{3a} t^2 \u }{\bc{\exp{a} t + b}^2}
  ,
  \qquad
  \g{\v} = \frac{\exp{2a} t^2 \v}{\bc{\exp{a} t + b}^2}
  + \frac{2\exp{2a}bt\, \u}{\bc{\exp{a} t + b}^3}
  .
\end{equation}
We note that the first equation in \eqref{eqn:b2} is not of the form
discussed in Section \ref{sec:genmeth}, as it contains the term $\u_t$.
Here we include terms in $\u_t$ to simplify the computations. The first order system
\begin{equation} \label{eqn:b2:altsys}
  \begin{split}
    t^2\v_{t} + 4t \v + 2 \u - \bc{2t\u + t^2\v}^{\frac12} & = 0,
    \\
    \v - \u_t & = 0,% \\
    %\u(0) = y_0, \quad \v(0) = y_1
    %,
  \end{split}
\end{equation}
would be equally valid. However, since the group action on
\eqref{eqn:b2:altsys} leads to significantly more complex expressions,
we do not consider this system here. In the spirit of \eqref{eqn:cg},
we introduce a standard finite element approximation by seeking for
$\U,\V \in \cpoly{q+1}\bc{I_n}$ such that
\begin{equation} \label{eqn:b2fe}
  \begin{split}
    \int_{I_n} \bc{
      t^2 \V_t + 4 t \U_t + 2 \U
      - \bc{ 2 t \U + t^2 \U_t}^{\frac12}
    } \phi \di{t}
    & = 0
    \qquad
    \forall \phi \in \dpoly{q}\bc{I_n},
    \\
    \int_{I_n} \bc{
      \U_t - \V
    } \psi \di{t}
    & = 0
    \qquad
    \forall \psi \in \dpoly{q}\bc{I_n}
    .
  \end{split}
\end{equation}
For any function $f(t) \in \dpoly{q}\bc{I_n}$ we observe that
$f\bc{g \cdot t}$ remains in $\dpoly{q}$ up to translation in the
time variable $t$.  Therefore, our space of test functions is
preserved up to temporal translations and we conclude that our space of
test functions is ``invariant'' under this group action as discussed
in Step 2.1 of Section \ref{sec:genmeth}. Acting on the finite element
approximation \eqref{eqn:b2fe} with \eqref{eqn:b2:sym} we find, after
simplification, the transformed functionals
\begin{equation} \label{eqn:b2:fe2}
  \begin{split}
    \int_{I_n} \exp{2a} \Bigg(
      t^2 \V_t
      + \frac{4 \exp{a} t^2 + 2b t}{\exp{a}t+b} \U_t
      + \frac{2bt \V}{\exp{a}t+b}
      + 2 \U
      - \bc{t^2 \U_t + 2t \U}^{\frac12}
    \Bigg) \phi \di{t}
    & = 0,
    \\
    \int_{I_n}
    \frac{\exp{3a} t}{\bc{\exp{a}t+b}^2}
    \bc{ \U_t - \V}
    \psi \di{t}
    & = 0
    .
  \end{split}
\end{equation}
We observe consistency of the first equation by equating
$\U_t$ and $\V$. However, the $\V$ term scales differently in $t$
from the second equation and so cannot be substituted.

To construct a moving frame, we choose the cross-section
\begin{equation} \label{eqn:b2:cs}
  \cs = \bw{\U = \textrm{sign}\bc{\U},\; \V = 0}
  .
\end{equation}
The resulting normalisation equations are
\begin{equation} \label{eqn:b2:norm}
  \frac{\exp{3a} t^2 \U }{\bc{\exp{a} t + b}^2}
  = \textrm{sign}\bc{\U}
  ,
  \qquad
  \frac{\exp{2a} t^2 \V}{\bc{\exp{a} t + b}^2}
  + \frac{2\exp{2a}bt \U}{\bc{\exp{a} t + b}^3}
  = 0
  .
\end{equation}
Solving for the group parameters $a,b$, under the assumptions that $\U
\neq 0$ and $\U \neq -\cfrac{t\V}{2}$, we obtain the moving frame
\begin{equation}
  a = \ln{\bigg(\frac{\abs{\U}}{\bc{\U + \frac12 t \V}^2}}\bigg)
  ,
  \qquad
  b = \frac{t^2 \U\V}{\bc{\U + \frac12 t\V}^3}
  .
\end{equation}
Substituting these normalised group parameters into \eqref{eqn:b2:fe2} yields the
invariant finite element approximation
\begin{equation} \label{eqn:b2fe2:invariant}
  \begin{split}
    \int_{I_n}
    \frac{\U}{\bc{\U+\frac12 t \V}^4} \bigg[
      \U \bc{
        t^2 \V_t + 4 t \U_t + 2 \U
        - \bc{t^2 \U_t + 2 t \U}^{\frac12}} \qquad \qquad & \\
      - t^2 \V \bc{
        \V - \U_t
      }\bigg] \phi \di{t}
      & = 0
    \quad
    \forall \phi \in \dpoly{q}\bc{I_n},
    \\
    \int_{I_n}
    \U^{-1} \bs{\U_t - \V} \psi \di{t}
    & = 0
    \quad
    \forall \psi \in \dpoly{q}\bc{I_n}
    ,
  \end{split}
\end{equation}
where $\U,\V \in \cpoly{q+1}\bc{I_n}$.
\end{example}

\begin{example}[{\bf A non-projectable action}] \label{example:noproject}
We now draw our attention to the first order ODE
\begin{equation} \label{eqn:noproject}
  \frac{y_t}{y-ty_t} - C = 0
  ,
\end{equation}
for a known fixed constant $C$, subject to the initial $y(0)=y_0$. This ODE
is invariant under the two-parameter non-projectable symmetry group action
\begin{equation} \label{eqn:noproject:sym}
  \g{t} = t + \alpha y,
  \qquad
  \g{y} = \exp{\beta} y
  ,
\end{equation}
where $\alpha,\beta\in \mathbb{R}$.
As the ODE is a first order system we do not need to introduce
auxiliary variables to propose a finite element
discretisation. However, for consistency with the remainder of this
work we shall write $y(t):=\u(t)$. A standard finite element
approximation for \eqref{eqn:noproject} is given by seeking for
$\U \in \cpoly{q+1}\bc{I_n}$ such that
\begin{equation} \label{eqn:noproject:fe}
  \int_{I_n} \bigg(
    \frac{\U_t}{\U-t\U_t} - C \bigg)
  \phi \di{t}
  = 0
  \qquad
  \forall \phi \in \dpoly{q}\bc{I_n}
  .
\end{equation}
Applying the prolonged group action 
\begin{equation}
  \g{t} = t + \alpha \U,
  \qquad
  \g{\U} = \exp{\beta} \U,
  \qquad
  \g{\U_t}
  = \frac{\exp{\beta} \U_t}{1 + \alpha \U_t}
  .
\end{equation}
 to the finite element approximation \eqref{eqn:noproject:fe} yields
\begin{equation} \label{eqn:noproject:action}
  \int_{I_n} \bigg(
    \frac{\U_t}{\U-t\U_t} - C\bigg)
  \bc{1 + \alpha\U_t} \g{\phi}
  \di{t} = 0
  .
\end{equation}
The induced action of the transformation $\g{t} = t+\alpha \U$ on the test functions $\phi$ will 
depend on the degree of the functions considered.
First, let us consider the case where the test functions are
naturally invariant, which occurs when $q=0$ and the test
functions are time independent. In this case $\g{\phi} = \phi$ and
we obtain the transformed finite element scheme
\begin{equation} \label{eqn:noproject:fe0}
  \int_{I_n} \bigg(
    \frac{\U_t}{\U-t\U_t} - C\bigg)
  \bc{1 + \alpha\U_t} \phi
  \di{t} = 0
  \qquad
  \forall \phi \in \dpoly{0}\bc{I_n}
  ,
\end{equation}
where $\U \in \dpoly{1}^C\bc{I_n}$.
To construct a moving frame, we choose the cross-section
\begin{equation} \label{eqn:nonproject:cs}
  \cs = \bw{t = 0,\; \U = \textrm{sign}\bc{\U}}
  .
\end{equation}
The corresponding normalisation equations are
\begin{equation}
  t + \alpha \U = 0,
  \qquad
  \exp{\beta} \U = \textrm{sign}\bc{\U}
  .
\end{equation}
Solving for the group parameters yields the moving frame
\begin{equation} \label{eqn:noproject:mf}
  \alpha = - \frac{t}{\U}
  ,
  \qquad
  \beta = - \ln\bc{\abs{\U}}
  .
\end{equation}
Substituting \eqref{eqn:noproject:mf} into \eqref{eqn:noproject:fe0} gives the 
invariant scheme where one seeks for $\U \in \cpoly{1}\bc{I_n}$ such that
\begin{equation} \label{eqn:noproject:invariant0}
  \int_{I_n}\bigg(
    \frac{\U_t - C \bc{\U-t\U_t}}{\U}
  \bigg) \phi \di{t} = 0
  \qquad \forall \phi \in \dpoly{0}\bc{I_n}
  .
\end{equation}
For $q>0$, the invariantisation of the finite element scheme is more
complicated, as the space of test functions is not invariant. In fact,
as the group action is not projectable there is no known conventional
choice of function space that is invariant under the symmetry group.
We must instead invariantise every basis function that spans our
function space, which is encapsulated into the general invariantisation
formula \eqref{eqn:invariantisation}. We note that after the
invariantisation of the basis functions we no longer expect them to
span a function space.  Nevertheless, implementing the
invariantisation procedure will produce a consistent
symmetry-preserving finite element method. As an illustrative
example, consider the case where $q=1$. Here the space of test functions
$\dpoly{1}\bc{I_n}$ may be given in terms of the Lagrange basis
functions
\begin{equation}\label{eqn:lagrangebasis}
  \dpoly{1}\bc{I_n}
  =
  \basis{\lagrange_1,\lagrange_2},
\end{equation}
where
\begin{equation} \label{eqn:linearlagrange}
  \lagrange_1\bc{t}
  =
  \frac{t-t_{n+1}}{t_n-t_{n+1}}
  \qquad\text{and}\qquad
  \lagrange_2\bc{t}
  =
  \frac{t-t_n}{t_{n+1}-t_n}
  .
\end{equation}
Letting
\[
\phi = a \lagrange_1(t) + b \lagrange_2(t),
\]
with $a,b\in \mathbb{R}$, the transformed finite element
method \eqref{eqn:noproject:action} becomes
\begin{equation}\label{eqn:q1:noproject:action}
  \int_{I_n} \bigg(
    \frac{\U_t}{\U-t\U_t} - C\bigg)
  \bc{1 + \alpha\U_t}
  \bc{a \mathcal{N}_1\bc{t,\U;\alpha}
    + b \mathcal{N}_2\bc{t,\U;\alpha} }
  \di{t} = 0
  ,
\end{equation}
where
\begin{equation}
  \begin{split}
    \mathcal{N}_1\bc{t,\U;\alpha}
    & =
    \frac{t - t_{n+1} + \alpha \bc{\U(t)-\U(t_{n+1})}}%
    {t_n-t_{n+1} + \alpha \bc{\U(t_n)-\U(t_{n+1})}},
    \\
    \mathcal{N}_2\bc{t,\U;\alpha}
    & =
    \frac{t - t_n + \alpha \bc{\U(t)-\U(t_n)}}%
    {t_{n+1}-t_n + \alpha \bc{\U(t_{n+1})-\U(t_n)}}
    ,
  \end{split}
\end{equation}
and $\U\in\cpoly{2}\bc{I_n}$.  Substituting the moving frame expressions \eqref{eqn:noproject:mf} into \eqref{eqn:q1:noproject:action} yields the following invariant finite element approximation: Seek
$\U\in\cpoly{2}\bc{I_n}$ such that for all $a,b\in \mathbb{R}$,
\begin{equation}
  \int_{I_n} \bigg(
    \frac{\U_t}{\U-t\U_t} - C\bigg)
  \bc{1 + \alpha\U_t}
  \bc{a \mathcal{M}_1\bc{t,\U}
    + b \mathcal{M}_2\bc{t,\U} }
  \di{t} = 0,
\end{equation}
where
\begin{equation}
  \begin{split}
    \mathcal{M}_1\bc{t,\U}
    & =
    \frac{t - t_{n+1} - \frac{t}{\U(t)} \bc{\U(t)-\U(t_{n+1})}}%
    {t_n-t_{n+1} - \frac{t}{\U(t)} \bc{\U(t_n)-\U(t_{n+1})}},
    \\
    \mathcal{M}_2\bc{t,\U}
    & =
    \frac{t - t_n - \frac{t}{\U(t)} \bc{\U(t)-\U(t_n)}}%
    {t_{n+1}-t_n - \frac{t}{\U(t)} \bc{\U(t_{n+1})-\U(t_n)}}
    .
  \end{split}
\end{equation}
We note that at the endpoints of the element, i.e.\ at $t_n$ and $t_{n+1}$, the
value of modified basis functions $\mathcal{M}_i$ is the same as with the Lagrange
basis functions.  However, we do not expect them to form a partition of
unity. This procedure may be replicated for arbitrarily high order
Lagrange basis functions, leading to arbitrarily high order invariant
finite element approximations.
\end{example}

\revise{
  \begin{remark}[{\bf The connection between Lie point symmetry preserving
    methods and other geometric numerical
    integrators}] \label{rem:geometricconnection}

  We observe that, in view of Remark \ref{rem:geometric}, under a
  $q+1$ point quadrature approximation the standard and invariant
  finite element approximations are equivalent, following the
  methodology outlined in the proof of \cite[Theorem
  3.1.9]{self:thesis}. This suggests that the perturbation we are
  making to our numerical method to preserve the Lie symmetry is
  small. In fact, within the context of Hamiltonian ODEs the standard
  finite element method is energy preserving, indicating that Lie
  point symmetry preserving integrators are ``close'' to energy
  preserving integrators. Furthermore, we note both finite element
  approximations under this quadrature yield a symplectic collocation
  method \cite[\S3.1]{self:thesis}, suggesting that Lie point symmetry
  preserving integrators are ``close'' to symplectic integrators. We
  conjecture that, within the field of geometric numerical
  integration, Lie point symmetry preserving schemes are competitive
  alternatives to energy preserving and symplectic
  methods. Significant theoretical work beyond the scope of this paper
  is required to verify this conjecture.

  \end{remark}
}

\FloatBarrier
%%%%%
\section{Numerical experiments} \label{sec:numerics}
%%%%%

In this section, we consider select numerical experiments showcasing
the impact of considering symmetry-preserving finite element
methods.
In particular, we assume that our element size
$\dt{n} = t_{n+1}-t_n=\dt{}$ is uniform. We note that there are no
theoretical requirements for considering a uniform element size. This
simplification is made solely to add clarity to the results of our
numerical experiments.

We measure all errors in the $L_2$ norm, which for a vectorial
solution is defined as
\begin{equation} \label{eqn:l2}
  \norm{L_2\bc{\bs{0,T}}}{\vec{U}-\vec{u}}
  :=
  \bigg(\sum_{i=0}^m \int_0^T \bc{U_i-u_i}^2 \di{t}\bigg)^{\frac12}
  ,
\end{equation}
where $\vec{U} = \bc{U_0,...,U_m}$ is the numerical solution and
$\vec{u}=\bc{u_0,...,u_m}$ is the corresponding exact solution. We
also consider the \revise{maximal error at the temporal nodes},
which we shall use to investigate nodal super-convergence. It is important to
remark that the maximal nodal error only induces a norm when $q=0$,
however, we shall not use it to investigate convergence rates in the
sequel instead relying on the $L_2$ norm \eqref{eqn:l2}.

We use the $L_2$ errors to compute an experimental order of
convergence (EOC), which is defined as follows.

\begin{definition}[{\bf Experimental order of convergence}] \label{def:eoc}
  Given two finite sequences $\{a_k\}_{k=0}^n$ and $\{b_k\}_{k=0}^n$, we define the experimental order
  of convergence (EOC) as
  \begin{equation} \label{eqn:eoc}
    EOC \left(a,b; k \right) = \frac{\log{\cfrac{a_{k+1}}{a_k}}}{\log{\cfrac{b_{k+1}}{b_k}}},\qquad k=0,\ldots,n-1.
  \end{equation}
  In the sequel $a_k$ will represent errors and $b_k$ will represent
  element sizes.
\end{definition}

In practice, our finite element approximation is solved through the
following steps:
\begin{enumerate}

\item The problem is linearised, allowing us to solve an underlying
  linear problem which converges to the solution of the nonlinear
  problem. In our numerical experiments we employ a Newton solver,
  which is solved up to a tolerance of $10^{-12}$.

\item Our finite element functions are decomposed into basis
  functions, for example
  $U(t) = U^1 \lagrange_1(t) + \cdots + U^{q+2} \lagrange_{q+2}(t)$,
  where $\lagrange_i(t)$ are the degree $q{+}1$ Lagrange basis
  functions over the element $I_n$ and $U^i$ are the values of the
  finite element function at the Lagrange points, otherwise known as
  the degrees of freedom. Additionally, decomposing the test functions
  into their basis functions we may assemble a linear system of
  equations $A \vec{U} = b$, where $\vec{U} = \bc{U^1,...,U^{q+2}}$,
  $A$ is a matrix that represents terms depending on both $U$ and the
  test functions, and $b$ is a vector involving terms depending solely
  on the test functions and time. We observe that for finite element
  formulations of the form \eqref{eqn:cg}, the linear system
  $A\vec{U}=b$ is square after the enforcement of initial data.

\item To solve the finite element approximation we must evaluate the
  integral, while it often may be computed exactly, in practice it is
  typically evaluated through a quadrature approximation. In the
  sequel we employ an order $16$ Gauss quadrature, the error of which we
  expect to be below machine precision when the quadrature is not
  exact.

\item Finally, we solve the linear system iteratively to obtain the
  finite element solution over one element.

\end{enumerate}

For more information on the practical implementation of finite element
methods see \cite[\S0]{BrennerScott:2007}. In the sequel we shall use
Firedrake, \cite{Firedrake:2017}, to conduct steps (2) and (3), and
utilise the NumPy \verb;linalg; routine \verb;solve;,
\cite{Numpy:2006}, for (4).

\revise{
  \begin{remark}[{\bf Computational costs}] \label{rem:quadcost}
    
    In step (3), we note that the quadrature approximation computed is of significantly higher order than the numerical
    method. While the quadrature approximation is of higher order, this does
    not significantly contribute to the computational
    complexity of the method. We also note that steps (1)--(3), up to the
    inclusion of initial data, may be conducted entirely independently
    of the particular time step under consideration. In fact, the
    linear system solved iteratively in step (4) is identical
    on each time step up to the enforcement of initial data. The
    linearisation of the method and cost of assembling basis functions
    and evaluating the integrals are an initial set-up cost of the
    method and may be conducted offline. In practice, we are only
    required to enforce different initial data on the linear system at
    each time step. Finally, we note that the Newton solver \emph{does} have a
    leading order computation cost as it is employed at every time
    step, which is typical of nonlinear numerical approximations.
    
  \end{remark}
}

% redefine \dt macro for tables
\renewcommand{\dt}{\ensuremath{\tau}}

\begin{example}[{\bf \revise{Working example}}] \label{numerics:painleve}

  Here we compare standard finite element approximation
  \eqref{eqn:painlevefe} of the \revise{working example}
  \eqref{eqn:painleve2} with the symmetry-preserving finite element
  scheme \eqref{eqn:painleve:invariant}, which is investigated as an
  illustrative example throughout Section \ref{sec:methodology}. We
  initialise both finite element schemes using
  \begin{equation} \label{eqn:num:ex21}
    \U(0) = 1, \qquad
    \V(0) = -1
    ,
  \end{equation}
  which approximates the exact solution
  \begin{equation} \label{eqn:num:ex22}
    \begin{split}
      \u(t)  = \exp{-t}, \qquad
      \v(t)  = - \exp{-t}
      ,
    \end{split}
  \end{equation}
  i.e., we wish to numerically simulate exponential decay, with this
  in mind, we consider a relatively short time simulation with $T=10$,
  as over significantly longer time the exact solution will, up to
  numerical precision, be zero.  Simulating the standard finite
  element scheme \eqref{eqn:painlevefe} for various polynomial degrees
  we obtain Table \ref{tab:ex20}. Similarly simulating the invariant
  finite element scheme \eqref{eqn:painleve:invariant} we obtain Table
  \ref{tab:ex21}.
  \begin{table}[h]
    \parbox{.45\linewidth}{%
      \captionsetup{width=\linewidth}
      \caption{The standard finite element approximation
        \eqref{eqn:painlevefe} where \eqref{eqn:num:ex21} and
        \eqref{eqn:num:ex22} hold with $T=10$. %We observe optimal
        % convergence in each polynomial degree.
        \label{tab:ex20}}
      \begin{tabular}{|c|c|c|c|c|}
\hline
 q & $\dt$    & Maximal  & $L_2$ error & EOC   \\
   &       	  &nodal error&			  & \\
  \hline
  & 1.56e-01 & 7.49e-04            & 1.70e-03    & - \\
 0 & 7.81e-02 & 1.87e-04            & 4.25e-04    & 2.00  \\
  & 3.91e-02 & 4.68e-05            & 1.06e-04    & 2.00  \\
  & 1.95e-02 & 1.17e-05            & 2.66e-05    & 2.00  \\
  \hline
  & 1.56e-01 & 3.04e-07            & 2.19e-05    & - \\
 1 & 7.81e-02 & 1.90e-08            & 2.74e-06    & 3.00  \\
  & 3.91e-02 & 1.19e-09            & 3.43e-07    & 3.00  \\
  & 1.95e-02 & 7.43e-11            & 4.28e-08    & 3.00  \\
  \hline
  & 1.56e-01 & 5.31e-11            & 1.58e-07    & - \\
 2 & 7.81e-02 & 8.30e-13            & 9.91e-09    & 4.00  \\
  & 3.91e-02 & 1.39e-14            & 6.20e-10    & 4.00  \\
  & 1.95e-02 & 4.75e-15            & 3.87e-11    & 4.00  \\
\hline
\end{tabular}
    }
    \hfill
    \parbox{.45\linewidth}{%
      \captionsetup{width=\linewidth}
      \caption{The invariant finite element approximation
        \eqref{eqn:painleve:invariant} where \eqref{eqn:num:ex21} and
        \eqref{eqn:num:ex22} hold with $T=10$. %We observe optimal
        % convergence in each polynomial degree, and that for all
        % simulations the approximation is exact at the end time.
        \label{tab:ex21}}
      \begin{tabular}{|c|c|c|c|c|}
\hline
 q & $\dt$    & Maximal  & $L_2$ error & EOC   \\
   &       	  &nodal error&			  & \\
  \hline
  & 1.56e-01 & 3.96e-16            & 2.23e-03    & - \\
 0 & 7.81e-02 & 2.84e-16            & 5.57e-04    & 2.00  \\
  & 3.91e-02 & 1.16e-15            & 1.39e-04    & 2.00  \\
  & 1.95e-02 & 7.77e-16            & 3.48e-05    & 2.00  \\
  \hline
  & 1.56e-01 & 5.83e-16            & 2.19e-05    & - \\
 1 & 7.81e-02 & 5.55e-16            & 2.74e-06    & 3.00  \\
  & 3.91e-02 & 7.77e-16            & 3.43e-07    & 3.00  \\
  & 1.95e-02 & 9.99e-16            & 4.28e-08    & 3.00  \\
  \hline
  & 1.56e-01 & 5.00e-16            & 1.58e-07    & - \\
 2 & 7.81e-02 & 1.17e-15            & 9.91e-09    & 4.00  \\
  & 3.91e-02 & 3.11e-15            & 6.20e-10    & 4.00  \\
  & 1.95e-02 & 4.77e-15            & 3.87e-11    & 4.00  \\
\hline
\end{tabular}
    }
  \end{table}
  First, we observe optimal convergence in each polynomial degree.  In
  addition, we notice that, while the $L_2$ errors in both the standard
  and invariant scheme are comparable,  the standard scheme
  has slightly smaller errors in the $L_2$ norm. Regardless, we
  observe that the invariant scheme is \emph{exact} at the nodes,
  indicating that by preserving the Lie point symmetries we exactly
  preserve the flow of the solution at the nodes. \revise{This phenomenon originates from the fact that the transformation group
    \eqref{eqn:painleve2:symU} is of the same form as the exact
    solution, and thereby allows for the exact preservation of the solution
    through exactly preserving the Lie point symmetries.}

  Now, consider the case where \eqref{eqn:painlevefe} models an
  exponential growth problem with initial conditions
  \begin{equation} \label{eqn:num:ex21bu}
    \U(0) = -1, \qquad
    \V(0) = -1
    ,
  \end{equation}
  and exact solution
  \begin{equation} \label{eqn:num:ex22bu}
    \begin{split}
      \u(t)  = \exp{t}, \qquad
      \v(t)  = \u(t)
      .
    \end{split}
  \end{equation}
  As both the invariant and non-invariant schemes are nonlinear, we
  note that this situation is difficult to simulate, as the
  exponential growth of the solution may prevent our nonlinear solver
  from converging for large time $t$.  This issue may be overcome by
  decreasing the step size $\dt{}$.  However, as $T$ increases
  linearly $\dt{}$ must decrease exponentially. Regardless, when $q=0$
  and $\dt{}=0.25$ we obtain Figure \ref{fig:ex2bu}. We observe that
  the invariant scheme is not only exact at the nodes, but that the
  error does not grow significantly over time. Conversely, we notice
  that the error of the standard scheme increases over time.
  \begin{figure}[h]
    \centering
    \caption{
      \revise{Absolute difference between the exact solution and
      the standard scheme \eqref{eqn:painlevefe} and the invariant
      scheme \eqref{eqn:painleve:invariant}, with $q=0$ and $\dt{}=0.25$.} 
      \label{fig:ex2bu}
    }
    \includegraphics[ width=0.75\textwidth,height=8.5cm]{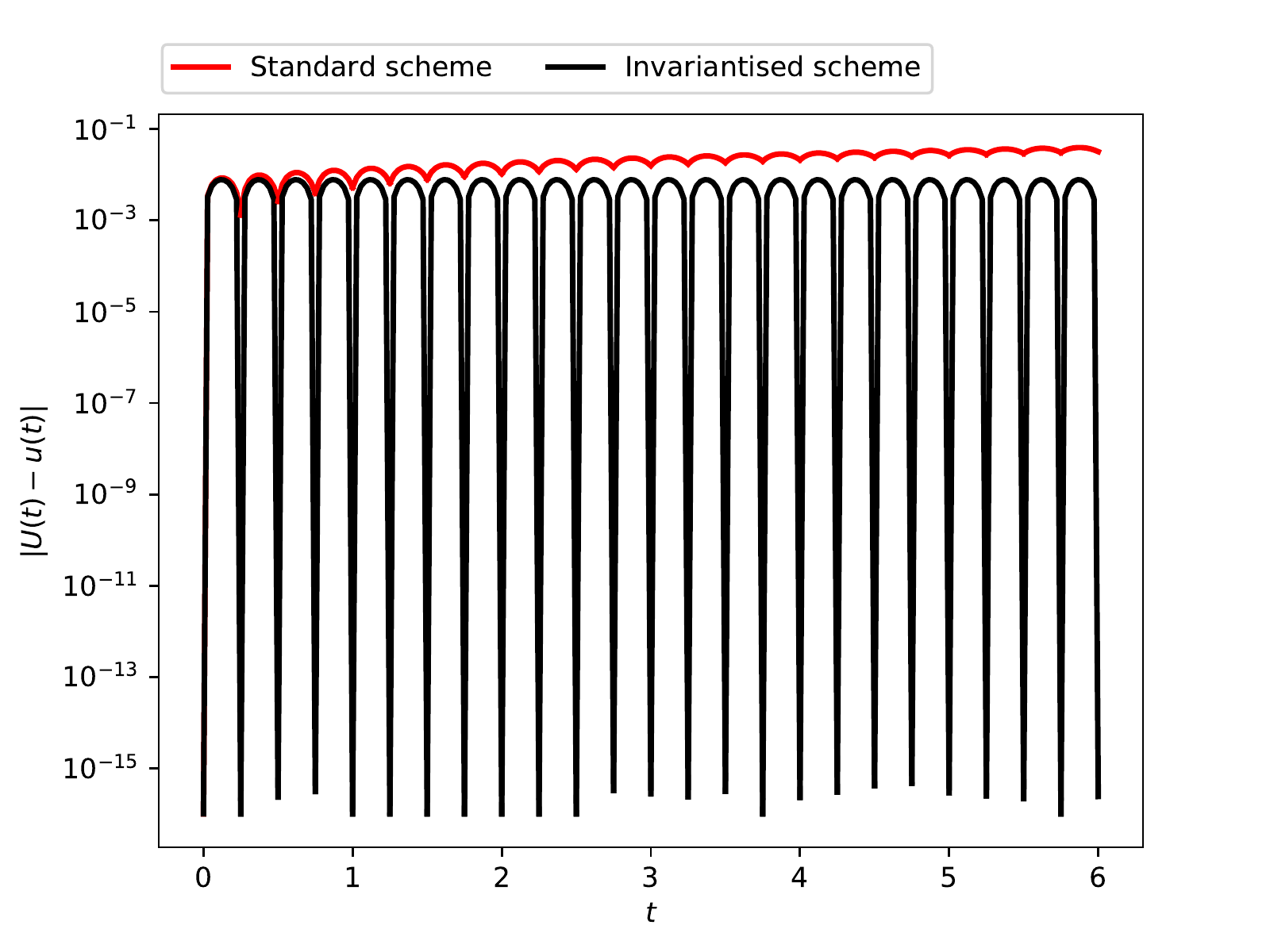}
  \end{figure}
  
\end{example}

\FloatBarrier
\begin{example}[{\bf Schwarzian differential equation}] \label{numerics:order3}

  We now compare the standard finite element scheme
  \eqref{eqn:order3:fe} against the invariantised approximation
  \eqref{eqn:order3:fe:invariant} for the third order ODE
  \eqref{eqn:order3} considered in Example \ref{example:order3}. In
  our numerical study we simulate the case where
  \begin{equation} \label{eqn:num:ex51}
    F(t) = 0
    ,
  \end{equation}
  subject to the initial conditions
  \begin{equation} \label{eqn:num:ex52}
    \U(0) = \W(0) = 1, \qquad
    \V(0) = -1
    ,
  \end{equation}
  up to the end time $T=1000$. Such numerical simulation approximates
  the exact solution
  \begin{equation} \label{eqn:num:ex53}
    \begin{split}
      \u(t) = \frac{4}{2+t} - 1 , \qquad
      \v(t)  = -\frac{4}{\bc{2+t}^2}, \qquad
      \w(t)  = \frac{8}{\bc{2+t}^3}
      .
    \end{split}
  \end{equation}
  The standard finite element approximation \eqref{eqn:order3:fe}
  leads to Table \ref{tab:ex50} while the invariant approximation
  \eqref{eqn:order3:fe:invariant} gives the results in Table
  \ref{tab:ex51}.
  
  \begin{table}[h]
    \parbox{.45\linewidth}{%
      \captionsetup{width=\linewidth}
      \caption{The standard finite element approximation
        \eqref{eqn:order3:fe} where \eqref{eqn:num:ex51},
        \eqref{eqn:num:ex52} and
        \eqref{eqn:num:ex53} hold with $T=1000$. %We observe optimal
        % convergence in each polynomial degree.
        \label{tab:ex50}}
      \begin{tabular}{|c|c|c|c|c|}
\hline
 q & $\dt$    & Maximal  & $L_2$ error & EOC   \\
   &       	  &nodal error&			  & \\
  \hline
  & 1.56e-01 & 1.45e-01       & 1.27e-01    & - \\
 0 & 7.81e-02 & 7.52e-02       & 3.17e-02    & 2.00  \\
  & 3.91e-02 & 3.83e-02       & 7.91e-03    & 2.00  \\
  & 1.95e-02 & 1.93e-02       & 1.98e-03    & 2.00  \\
  \hline
  & 1.56e-01 & 1.45e-01       & 7.79e-05    & - \\
 1 & 7.81e-02 & 7.52e-02       & 9.81e-06    & 2.99  \\
  & 3.91e-02 & 3.83e-02       & 1.23e-06    & 3.00  \\
  & 1.95e-02 & 1.93e-02       & 1.54e-07    & 3.00  \\
  \hline
  & 1.56e-01 & 1.45e-01       & 1.48e-06    & - \\
 2 & 7.81e-02 & 7.52e-02       & 9.38e-08    & 3.98  \\
  & 3.91e-02 & 3.83e-02       & 5.88e-09    & 4.00  \\
  & 1.95e-02 & 1.93e-02       & 3.68e-10    & 4.00  \\
\hline
\end{tabular}
    }
    \hfill
    \parbox{.45\linewidth}{%
      \captionsetup{width=\linewidth}
      \caption{The invariant finite element approximation
        \eqref{eqn:order3:fe:invariant} where \eqref{eqn:num:ex51},
        \eqref{eqn:num:ex52} and 
        \eqref{eqn:num:ex53} hold with $T=1000$. %We observe optimal
        % convergence in each polynomial degree.
        \label{tab:ex51}}
      \begin{tabular}{|c|c|c|c|c|}
\hline
 q & $\dt$    & Maximal  & $L_2$ error & EOC   \\
   &       	  &nodal error&			  & \\
  \hline
  & 1.56e-01 & 1.45e-01       & 3.60e-03    & - \\
 0 & 7.81e-02 & 7.52e-02       & 9.04e-04    & 1.99  \\
  & 3.91e-02 & 3.83e-02       & 2.26e-04    & 2.00  \\
  & 1.95e-02 & 1.93e-02       & 5.66e-05    & 2.00  \\
  \hline
  & 1.56e-01 & 1.45e-01       & 7.77e-05    & - \\
 1 & 7.81e-02 & 7.52e-02       & 9.81e-06    & 2.99  \\
  & 3.91e-02 & 3.83e-02       & 1.23e-06    & 3.00  \\
  & 1.95e-02 & 1.93e-02       & 1.54e-07    & 3.00  \\
  \hline
  & 1.56e-01 & 1.45e-01       & 1.48e-06    & - \\
 2 & 7.81e-02 & 7.52e-02       & 9.37e-08    & 3.98  \\
  & 3.91e-02 & 3.83e-02       & 5.88e-09    & 4.00  \\
  & 1.95e-02 & 1.93e-02       & 3.79e-10    & 3.95  \\
\hline
\end{tabular}
    }
  \end{table}
  For this problem, we observe that the maximal nodal error is
  comparable, however, the $L_2$ errors of the invariantised scheme
  are significantly smaller in the $q=0$ case highlighting a
  significant improvement through the invariantisation procedure. On the
  other hand, for
  higher polynomial degree, the $L_2$ errors are comparable for
  both schemes. Also, we note that in the $L_2$ norm, both schemes
  convergence optimally in each polynomial degree.
  \end{example}

\FloatBarrier
\begin{example}[{\bf A second order quasi-linear ODE}] \label{numerics:b}

We now consider the standard and invariant finite element schemes 
\eqref{eqn:b2fe} and \eqref{eqn:b2fe2:invariant}, which approximate the 
second order quasi-linear ODE \eqref{eqn:b}. We
enforce the initial data
\begin{equation} \label{eqn:num:ex41}
  \U(1) = 1, \qquad
  \V(1) = 2
  ,
\end{equation}
so that the exact solution is
\begin{equation} \label{eqn:num:ex42}
  \begin{split}
    \u(t) = \frac{t^3+9t^2+27t-25}{12t^2}, \qquad
    \v(t) = \frac{t^3-27t+50}{12t^3}
    ,
  \end{split}
\end{equation}
and simulate the solution over the domain $t \in [1,1000]$, where $T=1000$.
Simulating the standard finite element approximation \eqref{eqn:b2fe}
for various polynomial degree we obtain Table \ref{tab:ex40}, and for
the invariant approximation \eqref{eqn:b2fe2:invariant}, we obtain
Table \ref{tab:ex41}.
\begin{table}[h]
  \parbox{.45\linewidth}{%
  \captionsetup{width=\linewidth}
    \caption{The standard finite element approximation
      \eqref{eqn:b2fe} where \eqref{eqn:num:ex41} and
      \eqref{eqn:num:ex42} hold. %We observe optimal
      %convergence in each polynomial degree.
      \label{tab:ex40}}
    \begin{tabular}{|c|c|c|c|c|}
\hline
 q & $\dt$    & Maximal  & $L_2$ error & EOC   \\
   &       	  &nodal error&			  & \\
  \hline
  & 1.56e-01 & 2.41e-01       & 2.48e-02    & - \\
 0 & 7.81e-02 & 1.35e-01       & 6.30e-03    & 1.98  \\
  & 3.91e-02 & 7.25e-02       & 1.58e-03    & 1.99  \\
  & 1.95e-02 & 3.76e-02       & 3.96e-04    & 2.00  \\
  \hline
  & 1.56e-01 & 2.34e-01       & 1.22e-03    & - \\
 1 & 7.81e-02 & 1.34e-01       & 1.58e-04    & 2.94  \\
  & 3.91e-02 & 7.23e-02       & 2.00e-05    & 2.99  \\
  & 1.95e-02 & 3.75e-02       & 2.50e-06    & 3.00  \\
  \hline
  & 1.56e-01 & 2.34e-01       & 6.22e-05    & - \\
 2 & 7.81e-02 & 1.34e-01       & 4.11e-06    & 3.92  \\
  & 3.91e-02 & 7.23e-02       & 2.60e-07    & 3.98  \\
  & 1.95e-02 & 3.75e-02       & 1.64e-08    & 3.99  \\
\hline
\end{tabular}
  }
  \hfill
   \parbox{.45\linewidth}{%
   \captionsetup{width=\linewidth}
    \caption{The invariant finite element approximation
      \eqref{eqn:b2fe2:invariant} where \eqref{eqn:num:ex41} and
      \eqref{eqn:num:ex42} hold. %We observe optimal
      %convergence in each polynomial degree.
      \label{tab:ex41}}
    \begin{tabular}{|c|c|c|c|c|}
\hline
 q & $\dt$    & Maximal  & $L_2$ error & EOC   \\
   &       	  &nodal error&			  & \\
  \hline
 & 1.56e-01 & 2.43e-01       & 2.33e-02    & - \\
 0 & 7.81e-02 & 1.36e-01       & 6.09e-03    & 1.94  \\
 & 3.91e-02 & 7.25e-02       & 1.54e-03    & 1.98  \\
 & 1.95e-02 & 3.76e-02       & 3.87e-04    & 2.00  \\
  \hline
   & 1.56e-01 & 2.34e-01       & 1.26e-03    & - \\
 1 & 7.81e-02 & 1.34e-01       & 1.59e-04    & 2.99  \\
  & 3.91e-02 & 7.23e-02       & 2.00e-05    & 2.99  \\
  & 1.95e-02 & 3.75e-02       & 2.50e-06    & 3.00  \\
  \hline
  & 1.56e-01 & 2.34e-01       & 6.24e-05    & - \\
 2 & 7.81e-02 & 1.34e-01       & 4.10e-06    & 3.93  \\
  & 3.91e-02 & 7.23e-02       & 2.60e-07    & 3.98  \\
  & 1.95e-02 & 3.75e-02       & 1.67e-08    & 3.96  \\
\hline
\end{tabular}
  }
\end{table}
We observe that the errors for the standard and invariant scheme are
comparable, with a slight improvement in the errors in the case
$q=0$. These comparable results are a consequence of the exact
solution behaving linearly for large $t$, as linear behaviour may be
captured exactly by both of the schemes.
\end{example}

\FloatBarrier
\begin{example}[{\bf An ODE with non-projectable
  action}] \label{numerics:noproject}

The ODE
\begin{equation}
  \frac{\u_t}{\u-t\u_t} - C
  =
  0
  \qquad
  \u(0)=y_0
  ,
\end{equation}
has the simple solution given by
\begin{equation} \label{eqn:noproject:exact}
  \u(t) = y_0 \bc{Ct + 1}
  .
\end{equation}
Both the standard and invariant finite element approximations
experimentally achieve best approximability, i.e., they exactly
reconstruct polynomials of order $q{+}1$, as has been verified by the
authors numerically. This means that, in practice, even
for the case $q=0$ both finite element approximations are exact. As
such, making a comparison of their respective errors is an exercise in
futility. While the solution we are approximating is linear, both
finite element approximations are nonlinear.  That is to say that the
approximations must be linearised and iteratively solved. In our
experiments we utilised a Newton solver, which for sufficiently large
time steps may not converge. In Table \ref{tab:projectsolve} we have
tabulated a list of both standard and invariant schemes highlighting
at which time step size their respective Newton solvers fail to
converge.  We observe that the symmetry-preserving scheme allows for the problem to be solved with 
much larger time steps compared to its non-invariant counterpart.

\begin{table}[h]
  \caption{A table confirming whether the standard finite element approximation
    \eqref{eqn:noproject:fe} and the invariant approximation
    \eqref{eqn:noproject:fe0} may be successfully solved for various
    step sizes $\dt{}$ when approximating the exact solution
    \eqref{eqn:noproject:exact} with $C=1, y_0=0.5$.}
  \label{tab:projectsolve}
  \begin{tabular}{|c|c|c|}
  \hline
  $\dt$    & Standard scheme & Invariant scheme \\
  \hline
  0.390625 & \cmark          & \cmark \\
  0.78125  & \cmark          & \cmark \\
  1.5625   & \xmark          & \cmark \\
  3.125    & \xmark          & \cmark \\
  6.25     & \xmark          & \xmark \\
  \hline                               
\end{tabular}
\end{table}
\end{example}

\FloatBarrier

\revise{

\begin{example}[{\bf An illustrative example highlighting the benefits of
    preserving Lie point symmetries for a naive finite element
    discretisation}] \label{ex:naive}

  Excluding Example \ref{numerics:painleve}, where the benefits of
  invariantisation are clear, we note that the numerical experiments
  conducted thus far do not highlight
  overwhelming benefits of the invariantisation procedure. We
  conjecture that this is due to the natural preservation of geometric
  structures by the standard finite element method, as discussed in
  Remark \ref{rem:geometric}. Invariantising a less natural finite
  element approximation would indeed highlight the benefits of
  preserving Lie point symmetries, however, the authors feel this
  would make for an unfair comparison. For the sake of completeness,
  we now consider a poor finite element discretisation and show that
  its invariantisation recovers accurate long term simulations. We
  shall construct our naive approximation through the linearisation of
  an ODE model. Such an approximation may appear to be pathological, however, within the study of
  PDEs, linearisation is a common tool to obtain reduced models. This
  numerical example is motivated by the invariantisation of PDEs,
  which this work aims to inform.

  Consider the initial value problem
  \begin{equation} \label{eqn:vnl}
    \begin{gathered}
      y_{tt} = y^{-3}, \\
      y(0) = y_0, \quad y_t(0) = y_1
      ,
    \end{gathered}
  \end{equation}
  for some given constants $y_0, y_1$.  The ODE is invariant under the
  special linear group action
  \begin{equation} \label{eqn:vnl:sym}
    \g{t} = \frac{\alpha t + \beta}{\gamma t + \delta}
    ,
    \qquad
    \g{y} = \frac{y}{\gamma t + \delta}
    ,
    \qquad
    \alpha \delta - \beta \gamma = 1
    .
  \end{equation}
  Defining $\u=y$, we may rewrite \eqref{eqn:vnl} as the first
  order system
  \begin{equation} \label{eqn:vnl2}
    \begin{split}
      \v_t & = \u^{-3}, \\
      \u_t & = \v, \\
      \u(0) = y_0, &\qquad  \v(0) = y_1
      .
    \end{split}
  \end{equation}
  This system of ODEs is invariant under the extended group action
  \begin{equation} \label{eqn:vnl2:sym}
    \g{t} = \frac{\alpha t + \beta}{\gamma t + \delta}
    ,
    \qquad
    \g{\u} = \frac{\u}{\gamma t + \delta}
    ,
    \qquad
    \g{\v} = \bc{\gamma t + \delta}\v
    - \gamma \u
    .
  \end{equation}
  Now, instead of employing the finite element formulation \eqref{eqn:cg}, we now consider the 
  linearisation of \eqref{eqn:vnl}, which of course for a
  fundamentally nonlinear problem we expect will yield poor
  results. Thus, we formulate a finite element approximation
  of \eqref{eqn:vnl2} by seeking $\U,\V \in \cpoly{1}\bc{I_n}$ such that
  \begin{equation} \label{eqn:vnl2fe}
    \begin{split}
      \int_{I_n} \bc{
        \V_t - \U
      } \phi \di{t}
      & = 0
      \qquad
      \forall \phi \in \dpoly{0}\bc{I_n},
      \\
      \int_{I_n} \bc{
        \U_t - \V
      } \psi \di{t}
      & = 0
      \qquad
      \forall \psi \in \dpoly{0}\bc{I_n}
      .
    \end{split}
  \end{equation}
  This linearisation has been designed intentionally poorly, and is
  inconsistent with the solution of \eqref{eqn:vnl2}. Note that
  here we have fixed the degree of the finite element approximation for
  simplicity of exposition.  Higher order finite element approximations
  would yield similar results. In particular, as our test space is
  $\dpoly{0}\bc{I_n}$, i.e., the space of piece-wise constant
  functions, we note that the group action does not alter the test
  functions.
  
  Applying the group action \eqref{eqn:vnl2:sym} to the naive finite
  element scheme \eqref{eqn:vnl2fe} we obtain
  \begin{equation} \label{eqn:q0:nonlinearfe}
    \begin{split}
      \int_{I_n} \bs{
        \bc{\gamma t + \delta} \V_t
        - \frac{\U}{\bc{\gamma t + \delta}^3}
        + \gamma \bc{\V-\U_t}
      } \phi \di{t}
      & = 0, \\
      \int_{I_n} \bigg[
      \frac{\U_t - \V}{\gamma t + \delta}
      \bigg] \psi \di{t}
      & = 0
      ,
    \end{split}
  \end{equation}
  which shows that the finite element formulation is not invariant.  To obtain a symmetry-preserving formulation, we now implement the invariantisation process.  Assuming $\U\neq 0$, we choose the cross-section
  \begin{equation} \label{eqn:vnl2:crosssection}
    \cs = \bw{t = 0,\; \U = 1,\; \V = 0}
    ,
  \end{equation}
  and obtain the moving frame
  \begin{equation} \label{eqn:vnl:mf}
    \alpha = \U^{-1},
    \qquad
    \beta = -t \U^{-1},
    \qquad
    \gamma = \V,
    \qquad
    \delta = - t \V + \U
    .
  \end{equation}
  Substituting the group normalisations \eqref{eqn:vnl:mf} into \eqref{eqn:q0:nonlinearfe} we obtain the invariant finite element scheme consisting of seeking
  $\U,\V \in \cpoly{1}\bc{I_n}$ such that
  \begin{equation} \label{eqn:vnl2fe:invariant0}
    \begin{split}
      \int_{I_n} \bs{
        \V_t \U - \U^{-2}
        + \V\bc{\V-\U_t}
      } \phi \di{t}
      & = 0
      \qquad
      \forall \phi \in \dpoly{0}\bc{I_n},
      \\
      \int_{I_n} \U^{-1} \bc{
        \U_t - \V
      } \psi \di{t}
      & = 0
      \qquad
      \forall \psi \in \dpoly{0}\bc{I_n}
      .
    \end{split}
  \end{equation}
  We observe that by preserving the Lie point symmetries of \eqref{eqn:vnl}, we recover a consistent
  numerical scheme, even though
  consistency of the numerical scheme was destroyed by the initial
  naive discretisation. To illustrate the recovery of
  this long time structure we conduct a brief numerical
  experiment.  Setting $\dt{}=0.01$ and
  using the initial conditions
  \begin{equation} \label{eqn:num:ex31}
    \U(0) = 2^{\frac12}, \qquad
    \V(0) = 2^{-\frac12}
    ,
  \end{equation}
  corresponding to the exact solution
  \begin{equation} \label{eqn:num:ex32}
    \u(t) = \bc{t^2+2t+2}^{\frac12}, \qquad
    \v(t) = \bc{t+1}\bc{t^2+2t+2}^{-\frac12}
    ,
  \end{equation}
  we display the point-wise errors for the naive discretisation \eqref{eqn:vnl2fe} and the
  invariant discretisation \eqref{eqn:vnl2fe:invariant0} in Figure \ref{fig:ex3bu}.
    \begin{figure}[h]
    \centering
    \caption{\revise{ Absolute difference between the
        exact solution \eqref{eqn:num:ex32} and the
        naive discretisation \eqref{eqn:vnl2fe} and the
        invariant discretisation \eqref{eqn:vnl2fe:invariant0} with
        $\dt{}=0.01$.
        \label{fig:ex3bu}
      }
    }
    \includegraphics[ width=0.75\textwidth,height=8.5cm]{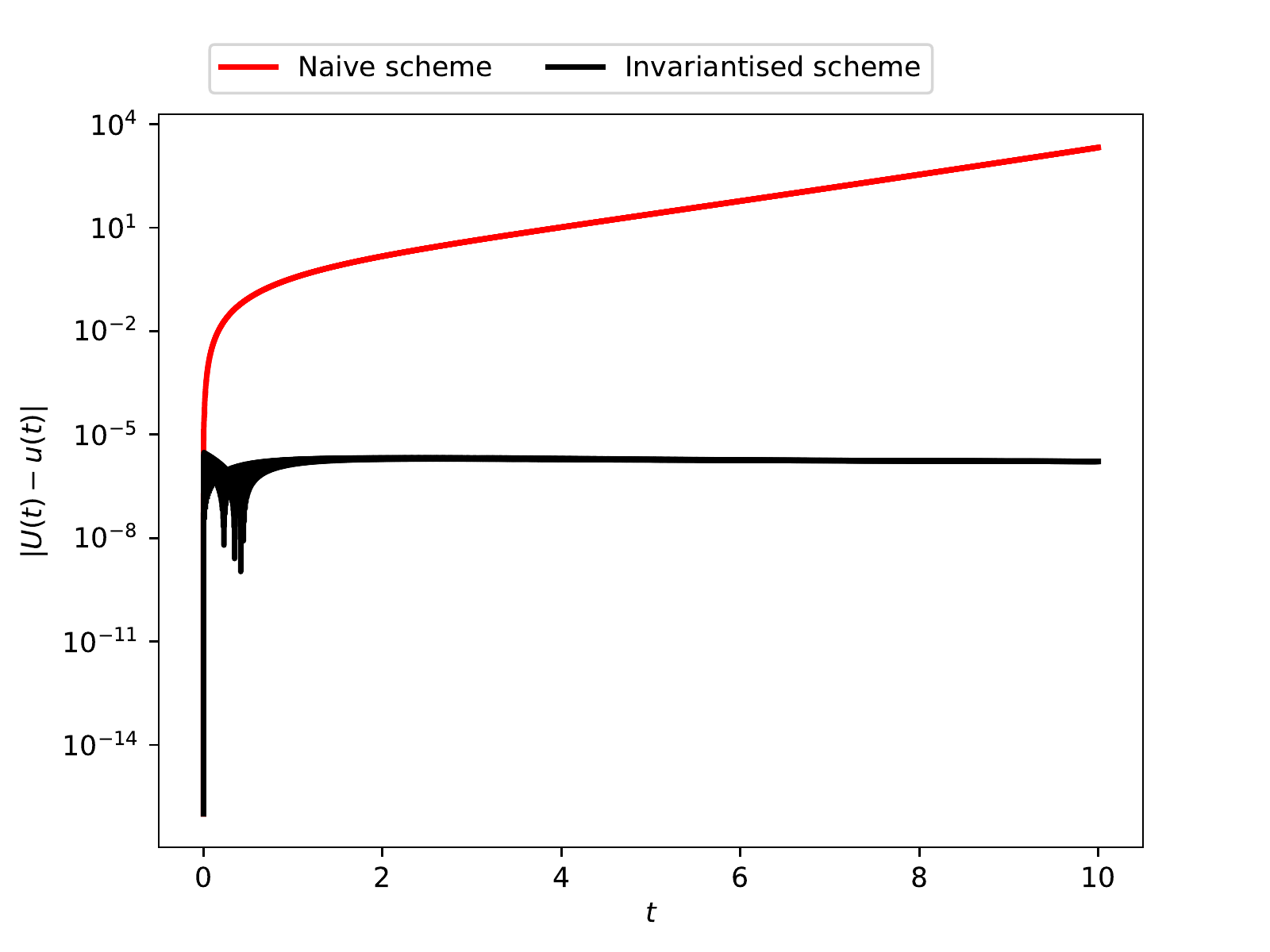}
  \end{figure}
  Similar to all previous examples in this section, we note that the
  invariant scheme \eqref{eqn:vnl2fe:invariant0} converges optimally,
  whereas the naive discretisation lacks consistency.
\end{example}
}
%\FloatBarrier

%%%%%
\section{\revise{A comparison of symmetry-preserving finite element methodologies}}\label{sec:comparison}
%\section{A comparison between the current methodology and an
%  existing methodology for finite element methods}
%%%%%

The preservation of Lie point symmetries for finite element methods
using the theory of moving frames has already been the subject of a
preliminary investigation in \cite{BihloValiquette:2018}. In this
final section, we highlight certain differences between the
methodology presented in this paper and the one considered in
\cite{BihloValiquette:2018}.

First, the equation used to implement the finite element method is
evidently different. In \cite{BihloValiquette:2018}, the finite
element method is applied to the original equation \eqref{eqn:ode}
while in the current work we consider the system of first order
equations \eqref{eqn:sys}. Obviously, the two formulations represent
the same equation, but the implementation of the finite element method
differs.  Indeed, when working with the original equation \eqref{eqn:ode},
it will be necessary to consider interpolating polynomials of degree
greater than $m$ (or interpolating functions with smoother boundary
conditions such as Hermite polynomials).  Therefore,
  as the order of the equation increases, the degree of the
  interpolating polynomials will increase and the level of
  computational difficulty will follow the same trend.  On the other
  hand, by recasting a differential equation as a system of first
  order equations, one can always work with low order interpolating
  polynomials.  Though, it is worth mentioning that as the order of the original ODE increases,
   the number of auxiliary variables $\vec{U}$ will also increase, which introduces it own set of
  computational challenges.

The second, and most important, distinction between the two approaches is in the interpretation of the finite element functions approximating the solution to the differential equation, and how the symmetry group acts on the interpolating functions.  In \cite{BihloValiquette:2018}, finite element methods are viewed in terms of their underlying difference discretisations.  As such, if 
\begin{equation}\label{eqn:Y}
Y(t) = \sum_n Y_n \phi_n(t)
\end{equation}
is an approximation of the exact solution $y(t)$ to the differential equation \eqref{eqn:ode}, the induced action on $Y(t)$ was defined as the combination of the product action on the coefficients $Y_n$, i.e.\ $\g{Y_n} = g\cdot Y_n$, together with the action on the basis functions $\g{\phi_n}$ introduced in \eqref{eqn:lift-phi} to give
\[
\g{Y} =\sum_n \g{Y_n}\g{\phi_n}.
\]
In this approach, time derivatives were approximated using finite differences and moving frames were constructed by normalising these discrete approximations, resulting in what are known as discrete moving frames, \cite{MariMansfield:2018}.  
For example, for the \revise{working example} \eqref{eqn:painleve} considered in Examples \ref{ex:painleve1}, \ref{ex:painleve2}, \ref{ex:painleve3} and \ref{ex:painleve4}, instead of considering the continuous 
cross-section \eqref{eqn:painleve:crosssection}, in \cite{BihloValiquette:2018} they introduce the discrete cross-section
\begin{equation}
  \cs = \bigg\{ Y_n = \textrm{sign}\bc{Y_n},\;
    \frac{Y_{n+1}-Y_{n-1}}{t_{n+1}-t_{n-1}} = 0 \bigg\}
  ,
\end{equation}
where the second term in the cross-section is the centred difference approximation 
of the first derivative $y_t$. We note that since
$Y(t)$ is linear and the second term in the cross-section spans two
elements, this cross-section is not equivalent to \eqref{eqn:painleve:crosssection}, although the two coincides in the continuous limit.  One of the issues of working with this approach is that when considering higher order
finite element schemes, one would need to specify cross-sections
utilising more accurate difference quotients to obtain high order invariant schemes.
On the other hand, in our
proposed framework this issue does not arise.

In contrast, in the present work we utilise
the continuous framework, and in particular, the continuously
differentiable nature of the finite element approximation when restricted
to the interior of a single element.  As such, the action on $\vec{U}(t)$ is defined to be the 
same continuous action as the one acting on the exact solution $u(t)$.  This significantly simplifies the
procedure and allows us to generate invariant schemes for arbitrary
degree approximations (assuming the action in $t$ is
linear). However, there are also certain drawbacks to our methodology. In
particular, as we do not evaluate the integrals of our approximation,
we may not remove time dependent factors multiplying every term in the
finite element approximation as in Examples 
\ref{ex:painleve2} and \ref{example:order3}. As such, we
find ourselves multiplying the entirety of our approximation by some
factor to ensure it is invariant under the integral. Such a problem
would not occur if considering the underlying finite difference
equation as in \cite{BihloValiquette:2018}.%, although the resultant
%methodology is highly dependent on the polynomial degree as we previously noted.

Another noticeable difference is in the group of transformations considered. In \cite{BihloValiquette:2018}, to preserve the form of the approximation \eqref{eqn:Y}, only projectable group actions were considered.  On the other hand, the current methodology applies to general Lie point symmetries.

Finally, we note that despite these differences, the two approaches
have shown that in certain cases the symmetry-preserving schemes
obtained can provide better long term numerical results than their
non-invariant counterparts. Understanding when this happens and the
varying benefits for different differential equations remains an open
question.

\FloatBarrier
\bibliographystyle{abbrv}
\bibliography{invariant}
\nocite{*}

\end{document}